\documentclass[twocolumn]{autart}    

\usepackage{graphics} 
\usepackage{epsfig} 
\usepackage{mathptmx} 
\usepackage{amsmath} 
\usepackage{amssymb}  
\usepackage{color}
\usepackage[noend]{algpseudocode}
\usepackage{url}
\usepackage{comment}
\usepackage{hologo}
\usepackage{cite}
\newtheorem{test}{TEST}

\newcommand{\supp}[1]{\operatorname{supp}(#1)}
\newcommand{\subsys}[2]{\ensuremath{(A, B_{(., {#1})}, C_{(#2, .)}, D_{(#2, #1)})}}

\begin{document}

\begin{frontmatter}

\title{ Securing State Estimation Under Sensor\\ and Actuator Attacks: Theory and Design } 

\thanks[footnoteinfo]{This paper was not presented at any IFAC 
meeting. Corresponding author M.~Showkatbakhsh Tel. +213 364-8655. }

\author[UCLA]{Mehrdad Showkatbakhsh}\ead{mehrdadsh@ucla.edu},    
\author[Maryland]{Yasser Shoukry}\ead{yshoukry@ece.umd.edu},               
\author[UCLA]{Suhas Diggavi}\ead{suhas@ee.ucla.edu},  
\author[UCLA]{Paulo Tabuada}\ead{tabuada@ucla.edu}  

\address[UCLA]{Electrical \& Computer Engineering Department, UCLA, Los Angeles, CA}  
\address[Maryland]{Electrical \& Computer Engineering, University of Maryland, College Park, MD}             


\begin{abstract}                          
This paper discusses the problem of estimating the state of a linear time invariant system when some of its sensors and actuators are compromised by an adversarial agent. In the model considered in this paper, the malicious agent attacks an input (output) by manipulating its value arbitrarily, i.e., we impose no constraints (statistical or otherwise) on how control commands (sensor measurements) are changed by the adversary. In the first part of this paper, we introduce the notion of sparse strong observability and we show that is a necessary and sufficient condition for correctly reconstructing the state despite the considered attacks. In the second half of this work, we propose an estimator to harness the complexity of this intrinsically combinatorial problem, by leveraging satisfiability modulo theory solving. Numerical simulations demonstrate the effectiveness and scalability of our estimator.
\end{abstract}

\end{frontmatter}

\section{Introduction}
Cyber-Physical Systems (CPS) are characterized by the tight interconnection of cyber and physical components. CPS are not only prone to actuator and sensor failures but also to adversarial attacks on the control and sensing modules. Security of CPS is no longer restricted to the cyber domain, and recent incidents such as the StuxNet malware \cite{stuxnet} and the security flaws reported on modern cars \cite{jeep, nissan} motivated the recent interest in security of CPS, (see for example, \cite{cardenas,sundaram1,par1_1,par1_2} and references therein). During the last decade, a number of security problems have been tackled by the control community, \emph{e.g.,} denial-of-service \cite{paulo_sonia,paulo_persis1,paulo_persis2,basar1}, replay attacks \cite{yilin3}, man-in-the-middle attacks \cite{paulo_roysmith}, false data injection \cite{yilin1}, etc. \\
This paper addresses the problem of state estimation when several {sensors \emph{and} actuators} are under attack. We broadly refer to state estimation in the adversarial environment as secure state estimation. Our attack model is quite general and we impose no constraints on the magnitude, statistical properties, or temporal characteristics of the signals manipulated by the adversary. \\
Secure state estimation has gained the attention of the control community over the past decade \cite{ACMsurvey18}. In one line of work, the problem of state estimation and control under sensor attacks is investigated and the authors derived necessary and sufficient conditions under which estimation and stabilization are possible \cite{fawzi}. Shoukry et. al. \cite{yasser1} further refined this condition and called it sparse observability. Chong et. al.  \cite{chong} found an equivalent condition for continuous-time systems and called it observability under attack. Nakahira et. a. \cite{yilin4} investigated a similar problem while considering the asymptotic correctness of state estimation. The authors relaxed the sparse observability condition to sparse detectability and showed it is a necessary and sufficient condition for asymptotic correctness. The noisy version of this problem has been investigated in the literature \cite{gupta1, gupta2, mo14, mo16, shaunak}. Mishra et. al. \cite{shaunak} derived the optimal solution for Gaussian noise. In this paper, we solve the more general problem of \emph{actuator and sensor} attacks that includes, as a special case, sensors attacks.\\
Under the sparse attack model in which an adversary can only target a bounded number of actuators and sensors, state estimation is intrinsically a combinatorial problem. Shoukry et. al.  \cite{yasser_SMT} proposed a novel secure state estimator using the Satisfiability Modulo Theory (SMT) paradigm, called \textsc{Imhotep-SMT}. The authors only considered attacks on sensors. In this paper we address the more general problem of sensor \emph{and} actuator attacks and build an SMT-based estimator that can correctly reconstruct the state under both types of attacks.\\
In another line of work, the problem of secure state estimation has been studied when the exact model of the system is not available \cite{yong, pajic}. Tiwari et. al. \cite{tiwari} proposed an online learning method by building so-called safety envelopes as it receives attack-free data to detect abnormality in the data when the system is prone to attacks. In \cite{showkatbakhsh1, showkatbakhsh2} the authors considered system identification under sensors attacks. In all of these works, the adversarial agent is restricted to only attacking sensors.\\
Pasqualetti et. al. \cite{fabio} investigated the problem of attack detection and identification. The authors related the undetectable and unidentifiable attacks to the zero-dynamics of the underlying system. The proposed attack identification mechanism consists of a number of fault-monitor filters that provide formal guarantees for the existence of the attack. The number of filters, however, grows exponentially with the number of attacked sensors/actuators, and therefore hinders scalability. In another work \cite{sandberg}, the authors investigated detectibility and identifiability of attacks in the presence of disturbances and the concept of security index is generalized to dynamical systems. The proposed method is inherently combinatorial and does not scale well with the number of attacked sensors and actuators. In this paper, by leveraging the SMT paradigm, we design a state estimator that scales well with the number of sensors and actuators.\\
Fault isolation and fault detection filters are classical control topics closely related to secure state estimation. The traditional fault tolerant filters can detect faults on actuators and sensors, however, they are not adequate for the purpose of security. Some of these filters assume a priori knowledge (statistical or temporal) of the fault signals \cite{faulttolerant}, an assumption that does not hold in the security framework. The classical fault detection filters \cite{jones} do not guarantee identification of all possible adversarial signals and zero-dynamics attacks remain stealthy. As an alternative approach, robustification has been used in order to estimate the state despite sparse attacks by either deploying Kalman filters or principle component analysis \cite{robust1, robust2}. The main drawback of these methods is the absence of formal guarantees for the correctness of the state. In contrast, the method proposed in this paper is guaranteed to construct the state correctly in spite of attacks on sensors \emph{and/or} actuators if the number of attacked components is below a specified threshold that depends on the system. In a recent work \cite{necmiye}, Harirchi et. al. proposed a novel fault detection approach using techniques from model invalidation. The authors pursued a worst-case scenario approach and therefore their framework is suitable for security. However, necessary and sufficient conditions for state estimation in a general adversarial setting were not investigated in \cite{necmiye}. In this paper, we precisely characterize the class of systems, by providing necessary and sufficient conditions, for which state reconstruction is possible despite sensor and/or actuator attacks.\\
The contributions of this paper can be summarized as follows:
\begin{itemize}
\item We introduce the notion of sparse strong observability by drawing inspiration from sparse observability \cite{fawzi,yasser1} and the classical notion of strong observablity \cite{strong_obs}. We show this is the relevant property when the adversarial agent not only compromises sensor measurements but can also attack inputs. 
\item We develop an observer by leveraging the SMT approach to harness the exponential complexity of the problem. Our observer consists of two blocks interacting iteratively until the true state is found (see Section \ref{section:secure_observer} for a detailed explanation of the observer's architecture). 
\item We propose two methods to further decrease the running time of the proposed algorithm by reducing the number of iterations of the observer. The first method exploits heuristics that can be efficiently computed at each iteration. The second method is inspired by the \textsc{QuickXplain} algorithm  \cite{quickxplain} that efficiently finds an irreducibly inconsistent set (see Section \ref{section:secure_observer} for a detailed discussion on the aforementioned methods). We demonstrate the scalability of our proposed observer by several numerical simulations.
\end{itemize}
A preliminary version of some of the results in this paper were presented in \cite{showkatbakhsh3} where we introduced the notion of sparse strong observability and drew the connection to secure state estimation. However, the formal proofs were not provided due to space limitations. Furthermore, we propose a new observer that outperforms the observer introduced in \cite{showkatbakhsh3}. This paper is organized as follows. Section \ref{section:prob_def} introduces notation followed by the attack model and the precise problem formulation. In Section \ref{section:conditions}, we introduce the notion of sparse strong observability and relate this notion to the problem of state reconstruction when some of the inputs and outputs are under adversarial attacks. This section concludes with the main theoretical contribution of this paper that is Theorem \ref{thm:main}. Section \ref{section:secure_observer} is devoted to designing an observer by exploiting the SMT paradigm. Section \ref{section:simulation_results} provides the simulation results followed by Section \ref{section:conc} that concludes the paper.

\section{Problem Definition}
\label{section:prob_def}
\subsection{Notation}
We denote the sets of real, natural and binary numbers by $\mathbb{R}$, $\mathbb{N}$ and $\mathbb{B}$. We represent vectors and real numbers by lowercase letters, such as $u$, $x$, $y$, and matrices with capital letters, such as $A$. Given a vector $x \in \mathbb{R}^n$ and a set \mbox{$O \subseteq \{1, \hdots, n \} $}, we use $x \vert_{O}$ to denote the vector obtained from $x$ by removing all elements except those indexed by the set $O$. Similarly, for a matrix $C \in \mathbb{R}^{n_1 \times n_2}$ we use $C\vert_{(O_1 ,O_2 )}$ to denote the matrix obtained from $C$ by eliminating all rows and columns except the ones indexed by $O_1$ and $O_2$, respectively, where $O_i \subseteq \{1, \hdots, n_i\}$ with $n_i \in \mathbb{N}$ for $i \in \{1, 2\}$. In order to simplify the notation, we use $ C\vert_{(. ,O_2)} := C\vert_{(\{1, \hdots,n_1 \}, O_2 )}$ and $ C\vert_{(O_1, .)} := C\vert_{(O_1, \{1, \hdots,n_2 \} )}$. We denote the complement of $O$ by \mbox{$\overline{O} := \{1, \hdots, n\} \setminus O$}. We use the notation $\{x(t)\}_{t=0}^{T-1}$ to denote the sequence \mbox{$x(0), \hdots, x(T-1)$}, and we drop the sub(super)scripts whenever it is clear from the context.

A Linear Time Invariant (LTI) system is described by the following equations:
\begin{align}
x(t+1) &= Ax(t) + Bu(t), \nonumber \\
y(t) &= Cx(t) + Du(t), \label{generic_sys}
\end{align} where $u(t) \in \mathbb{R}^m$, $x(t) \in \mathbb{R}^n$ and $y(t) \in \mathbb{R}^p$ are the input, state and output variables, respectively, $t \in \mathbb{N} \cup \{ 0 \}$ denotes time, and $A$, $B$, $C$ and $D$ are system matrices with appropriate dimensions. We use $(A, B,C ,D)$ to denote the system described by \eqref{generic_sys}. The order of an LTI system is defined as the dimension of its state space. A trajectory of the system consists of an input sequence with its corresponding output sequence. For an LTI system, 
\begin{align}
\mathcal{O}_{(A,C)} &:= \begin{bmatrix}
C^T & A^TC^T & \hdots & (A^T)^{n-1} C^T
\end{bmatrix}^T, \\
\mathcal{N}_{(A,B,C,D)} &:= \begin{bmatrix}
D & 0 &  \hdots & 0 \\
CB & D &  \hdots & 0\\
\vdots & & \ddots  &\\
CA^{n-2}B & CA^{n-3}B & \hdots & D \\
\end{bmatrix},
\end{align} are the \emph{observability} and \emph{invertibility} matrices, respectively, where $n$ is the order of the underlying system. In this paper, we often work with subsets of inputs and outputs. For a subset of outputs $\Gamma_{y} \subseteq \{1, \hdots, p\}$, we use the notation \mbox{$\mathcal{O}_{\Gamma_{y}} := \mathcal{O}_{(A,C\vert_{(\Gamma_{y}, . )})}$} to denote the observability matrix of outputs in the set ${\Gamma_{y}}$. For a set of inputs ${\Gamma_{u}} \subseteq \{1, \hdots, m\}$, we use the notation $\mathcal{N}_{{\Gamma_{u}} \to {\Gamma_{y}}}$ to denote $\mathcal{N}_{\subsys{{\Gamma_{u}}}{{\Gamma_{y}}}}$. 
For $x \in \mathbb{R}^n$, we define its support set as the set of indices of its non-zero components, denoted by $\supp{x}$. Similarity we define the support of the sequence $\{ x(t) \}$ as $\supp{\{x(t)\} } := \cap_{t} \supp{x(t)} $.
The observer proposed in this paper uses batches of inputs and outputs in order to reconstruct the state. We reserve capital bold letters to denote these batches,
\begin{align}
\bold{Y}^{\tau}(t) &:= \begin{bmatrix}  y(t-\tau+1)^T & \hdots & y(t)^T \end{bmatrix}^T, \\
\bold{U}^{\tau}(t) &:= \begin{bmatrix}  u(t-\tau+1)^T & \hdots & u(t)^T \end{bmatrix}^T,
\end{align} where $\tau \leq n$. Whenever $\tau$ is the order of the underlying system, we may drop the superscript for ease of notation. For a subset of outputs (inputs), denoted by $\Gamma_y \subseteq \{1, \hdots, p\}$ ($\Gamma_u \subseteq \{1, \hdots, m\}$), we use the notation $\bold{Y}^{\tau}\vert_{\Gamma_y}(t)$ ($\bold{U}^{\tau}\vert_{\Gamma_u}(t)$) for the batches of length $\tau$ that only consists of outputs (inputs) in the set $\Gamma_y$ ($\Gamma_u$).
For a vector $x \in \mathbb{R}^n$, we denote a generic norm, $l_2$-norm and $l_1$-norm of $x$ by $\|x \|$, $\|x \|_2$ and $\|x \|_1$. 

\subsection{System and Attack model}

This work is concerned with the problem of state reconstruction of LTI systems. We consider the scenario in which sensors and actuators are both prone to adversarial attacks. The ultimate goal is to reconstruct the state despite these attacks. In this part, we define the attack model and conclude this section with the precise problem statement.
The system $S$, is described by the following equations:
\begin{align}
x(t+1) = Ax(t) + Bu_S(t) , \nonumber \\
y_S(t) = Cx(t) + Du_S(t). \label{sys}
\end{align} Without loss of generality we assume $\begin{bmatrix} B^T & D^T \end{bmatrix}^T$ to be of full column rank.

\begin{figure}[h]
\centering
\includegraphics[scale=0.45]{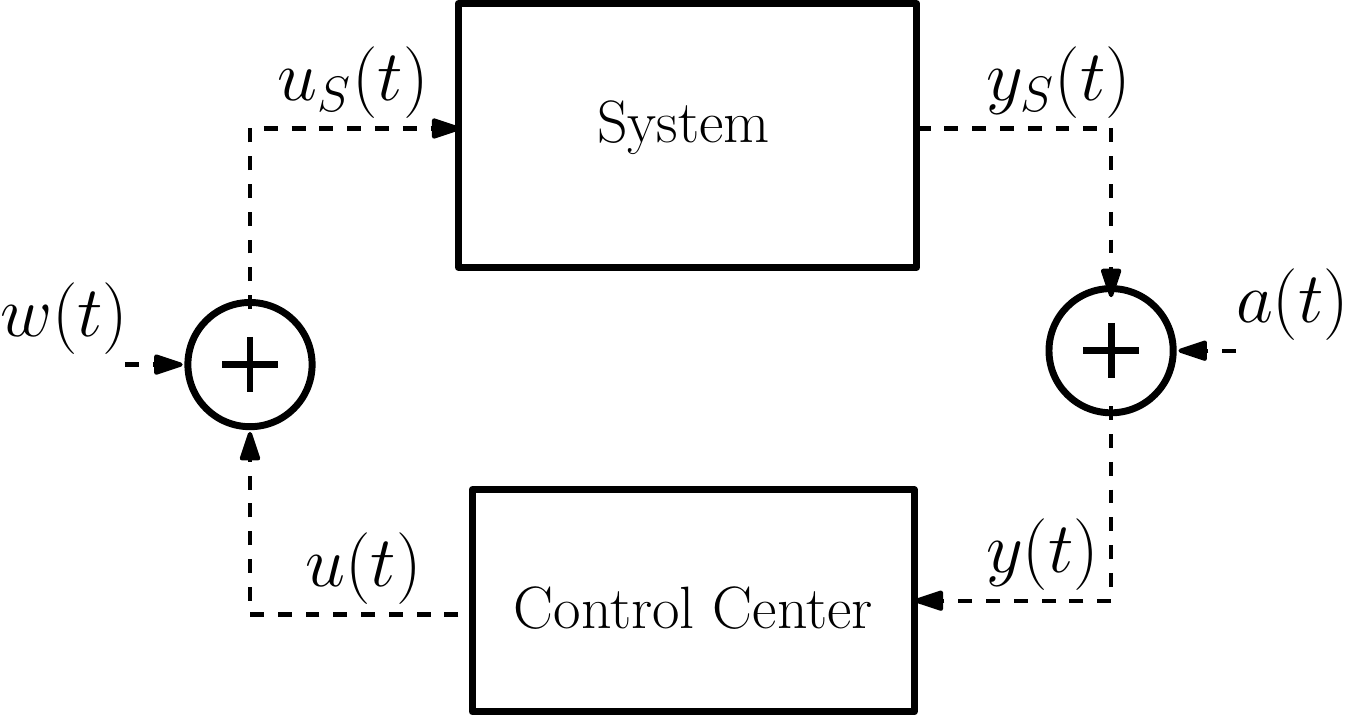}
\caption{The generic attack model considered in this paper.}
\label{fig:sys}
\end{figure}

Each actuator (sensor) corresponds to one input (output) and we use input (output) instead of actuator (sensor) in the rest of this paper. In this set up the adversary can attack both inputs and outputs. We model these attacks by additive terms and by imposing a sparsity constraint on them,
\begin{align}
\begin{cases}
u_S(t) &= u(t) + w(t),  \\
y(t) &= y_S(t) + a(t), \label{model}
\end{cases}
\end{align} where $u(t) \in \mathbb{R}^m$ and $y(t) \in \mathbb{R}^p$ are the controller-designed input and the observed output, respectively, and $w(t) \in \mathbb{R}^m$ and $a(t) \in \mathbb{R}^p$ are signals injected by the malicious agent. In the rest of this paper, we refer to these signals $\left( w \left( t \right), a(t) \right)$ as the attack of the adversarial agent. We use the subscript $S$ for signals that directly come from/to the system. The controller can only observe $y(t)$ and compute the input $u(t)$. This generic attack model is depicted in Figure \ref{fig:sys}.\\
When the adversary attacks an input (output) it can change its value to any arbitrary number without explicitly revealing its presence. The only limitation that we impose on the power of the malicious agent is the maximal number of inputs and outputs that can be attacked. 
\begin{assum}[Bound on the number of attacks]
The number of inputs and outputs under attack are bounded by $r$ and $s$, respectively.
\label{assump}
\end{assum}
Therefore, the malicious agent can attack a subset of inputs and outputs denoted by $\Gamma_u \subseteq \{1, \hdots, m \}$ and $\overline{\Gamma}_y \subseteq \{1, \hdots, p\}$,\footnote{For ease of exposition, we use $\Gamma_u$ to denote under-attack inputs while using $\Gamma_y$ for the set of attack-free outputs, i.e., the set of under-attack outputs is represented by $\overline{\Gamma}_y := |\{1, \hdots, p\} \setminus \Gamma_y$ in this paper. } respectively, with $|\Gamma_u| \leq r$ and $|\overline{\Gamma}_y| \leq s$, such that $\supp{\{w(t)\}} \subseteq \Gamma_u$ and $\supp{\{a(t)\}} \subseteq \overline{\Gamma}_y$. Note that these sets are not known to the controller and only upper bounds on their cardinality are given. Once the adversary chooses these sets, inputs and outputs outside these sets remain attack-free. This assumption is realistic when the time it takes for the adversarial agent to attack new inputs and outputs is large compared to the time scale of the system.\\
We now precisely define the main problem we tackle in this paper.
\begin{prob}[Secure state estimation]
For the linear system defined by \eqref{sys} under the attack model defined by \eqref{model}, what are necessary and sufficient conditions under which the state of the compromised system \eqref{sys} can be reconstructed with bounded delay? \label{problem:theory}
\end{prob} 
It is well-known that the secure state estimation problem, when only outputs are under adversarial attacks, is combinatorial and belongs to the class of \textit{NP-hard} problems \cite{yasser_SMT, fabio}. Therefore we are motivated to design an observer that harness the complexity of this problem.
\begin{prob}[Secure observer design]
Assumming conditions in Problem \ref{problem:theory} are satsified, how can we design an observer that reconstructs the state of the compromised system? \label{problem:design}
\end{prob}
\section{Conditions for Secure State Estimation }
\label{section:conditions}
In this section, we solve Problem \ref{problem:theory}, i.e., we provide conditions on the system described by \eqref{sys} under which state reconstruction (with bounded delay) is possible. We first develop the notion of sparse strong observability. This section concludes with Theorem \ref{thm:main} that relates this notion to the solution of Problem \ref{problem:theory}.\\
In the absence of attacks, the problem of estimating the state of a system while some of the inputs are unknown has been studied and the notion of strong observability was introduced in the literature \cite{strong_obs}. For strongly observable systems, it is possible to estimate the state of the system without the knowledge of inputs. The following definition formalizes this concept.
\begin{defn}[Strong observability]
An LTI system is called strongly observable if for any initial state $x(0) \in \mathbb{R}^n$ and any input sequence $\{ u(t) \in \mathbb{R}^m \}_{t = 0}^{\infty}$ there exists an integer $\tau \in \mathbb{N} \cup \{ 0 \}$ such that $x(0)$ can be uniquely recovered from $\{ y(t) \}_{ t = 0 }^{ \tau  }$.
\end{defn}
Note that $\tau$ is always upper-bounded by the order of the system. Linearity implies the following lemma.
\begin{lem}
An LTI system is strongly observable if and only if $y(t) = 0$ $\forall t \in \mathbb{N} \cup \{ 0 \}$ implies that $x(0) = 0$.\label{cor:strong_observable0}
\label{lemma:pre}
\end{lem}
\begin{pf}
Please refer to Appendix.
\end{pf}
It is straightforward to conclude the following corollary.
\begin{cor}
An LTI system is not strongly observable if and only if there exist a non-zero intial state and an input sequence such that $y(t) = 0$ for $t \in \mathbb{N} \cup \{ 0 \}$. \label{cor:strong_observable1}
\end{cor}
\begin{pf}
Follows directly from Lemma \ref{lemma:pre}.
\end{pf}
It is well-understood that when the adversary is restricted to attacking outputs, state reconstruction is possible only if there is enough redundancy in the outputs of the system. This redundancy can be stated in terms of observability of the system while removing a number of outputs. This property has been formalized in \cite{fawzi} and is called sparse observability \cite{yasser1}. By analogy with sparse observability, we define the notion of $(r,s)$-sparse strong observability as follows:
\begin{defn}[$(r,s)$-sparse strong observability]
An LTI system $(A, B, C, D)$ with $m$ inputs and $p$ outputs is $(r,s)$-sparse strongly observable if for any $\Gamma_{u} \subseteq \{1, \hdots, m\}$ and ${\Gamma}_{y} \subseteq \{1, \hdots, p\}$ with $| \Gamma_{u} | \leq r$ and $| {\Gamma}_{y} | \geq p-s$, the system $\subsys{\Gamma_{u}}{\Gamma_{y}}$ is strongly observable. \label{def:sparp_strong}
\end{defn}

Note that in Definition \ref{def:sparp_strong}, the value of $r$ and $s$ are upper bounded by the number of inputs and outputs, respectively. This modified notion of strong observability is the key for formalizing redundancy across inputs and outputs. We show that a necessary and sufficient condition for secure state estimation can be stated using this property. Note that $(0, s)$-sparse strong observability is equivalent to the notion of $s$-sparse observability that was introduced before in the literature \cite{fawzi,yasser1, shaunak}. The following theorem is the main theoretical result in this paper.

\begin{thm}
Let the number of attacked inputs and outputs be bounded by $r$ and $s$, respectively. Under the attack model \eqref{model}, the state can be reconstructed (possibly with delay) if and only if the underlying system is $(2r, 2s)$-sparse strongly observable. \label{thm:main}
\end{thm}

\begin{rem}
It is worth mentioning that the maximum number of attacked outputs, $s$, cannot be greater than $\left \lfloor{ \frac{p}{2} }\right \rfloor$ and it is an inherent limitation of LTI systems with $p$ outputs \cite{fawzi}. However the maximum number of attacked inputs is not inherently restricted by $\left \lfloor{ \frac{m}{2} }\right \rfloor$ and can take values up to $m$, depending on the specific system under the consideration.
\end{rem}
\begin{rem}
Pasqualetti et. al. \cite{fabio} addressed the problem of attack detection and identification in the presence of adversarial inputs and outputs for continuous-time LTI systems. They showed that attack identification is possible if and if  for any $\Gamma_{u} \subseteq \{1, \hdots, m\}$ and ${\Gamma}_{y} \subseteq \{1, \hdots, p\}$ with $| \Gamma_{u} | \leq 2r$ and $| {\Gamma}_{y} | \geq p-2s$, the system $\subsys{\Gamma_{u}}{\Gamma_{y}}$ does not have any invariant zeros.\\
It is clear that from the state and the dynamics of the system, the attack can be identified, therefore the attack identification comes free with the solution to the secure estimation problem. Strongly observable LTI systems do not have any invariant zeros (see, for example Theorem 1.8 in \cite{strong_obs}). Therefore this theorem shows that under this sparse-attack model, the conditions for identifying the attack also enable one to reconstruct the state, i.e., characterizations of attack identifiability and secure state estimation are equivalent for LTI systems. Putting these together, secure state estimation also comes with the solution to the attack identification problem. However, we provide a direct proof that does not require this machinery.
\end{rem}
\begin{pf}
First we show that $(2r, 2s)$-sparse strong observability is a sufficient condition for correctly estimating the state. For the sake of the contradiction, assume that the state cannot be reconstructed, i.e., there exist two different (initial) states, denoted by $x^{(1)}$ and $x^{(2)}$, that cannot be distinguished under this attack model. More precisely, there exist two attack strategies that will lead to the same exact (observed) trajectories. We reserve superscripts $.^{(1)}$ and $.^{(2)}$ for variables across those scenarios. Let us denote the adversarial additive terms by $\{w^{(1)}(t) \}, \{ a^{(1)}(t) \}$ and $\{ w^{(2)}(t) \}, \{ a^{(2)}(t) \}$. We represent the corresponding inputs and outputs of the system by $\{ u^{(1)}_{S}(t) \}, \{ y^{(1)}_{S}(t) \}$ and $\{ u^{(2)}_{S}(t) \}, \{ y^{(2)}_{S}(t) \}$, and the common (corrupted) measured output and the controller input sequences are denoted by $\{y(t)\}$ and $\{ u(t) \}$, respectively.\\
By the assumption of the attack model \eqref{model}, there exist $\Gamma_{u}^{(i)}, \overline{\Gamma}_{y}^{(i)}$ for $i \in \{1, 2\}$ with bounded cardinality such that
\begin{align}
\supp{ \{ w^{(i)}(t) \} } \subseteq \Gamma_{u}^{(i)}, \supp{ \{ a^{(i)}(t) \} } \subseteq \overline{\Gamma}_{y}^{(i)},
\end{align} for $i \in \{1, 2\}$. Note that
\begin{align}
 \begin{cases}  u^{(1)}_S(t)  = u(t) +  w^{(1)}(t) \\
									  u^{(2)}_S(t)  = u(t) +  w^{(2)}(t) \end{cases},
\end{align} where $u(t)$ is the controller designed input. Therefore 
\begin{align}
\supp{ \{ u^{(1)}_S(t) - u^{(2)}_S(t) \} } &= \supp{ \{ w^{(1)}(t) - w^{(2)}(t) \} }  \nonumber \\ 
& \subseteq \Gamma^{(1)}_{u} \cup \Gamma^{(2)}_{u}. 
\end{align}
Similarly, it is straightforward to conclude that {$\supp{ \{ y^{(1)}_S(t) - y^{(2)}_{S}(t) \} } \subseteq  \overline{\Gamma}^{(1)}_{y} \cup \overline{\Gamma}^{(2)}_{y}$}. We are ready to reach the contradiction. The underlying system is LTI, thus the input sequence $\{u^{(1)}_{S}(t) - u^{(2)}_{S}(t) \}$ with the initial state $x^{(1)}-x^{(2)}$ generates the output sequence $\{y^{(1)}_{S}(t) - y^{(2)}_{S}(t) \}$. The underlying system is $(2r, 2s)$-sparse strongly observable so the sub-system $\subsys{\Gamma_{u}}{\Gamma_{y}}$ is strongly observable for any $|\Gamma_{u} | = 2r$ and $|{\Gamma}_{y} | = p - 2s$. Let us choose $\Gamma_{u}$ and ${\Gamma}_{y}$ as any set of $2r$ inputs and $p -2s$ outputs such that,
\begin{align}
\Gamma^{(1)}_{u} \cup \Gamma^{(2)}_{u} \subseteq  \Gamma_{u}, \quad  \Gamma_{y} \subseteq {\Gamma}^{(1)}_{y} \cap {\Gamma}^{(2)}_{y}.
\end{align}
Note that $\{y^{(1)}_{S}(t)\vert_{ {\Gamma}_{y} } - y^{(2)}_{S}(t)\vert_{ {\Gamma}_{y} } \}$ is a zero sequence, hence by Lemma \ref{cor:strong_observable0} we conclude that the corresponding initial state ($x^{(1)} - x^{(2)}$) is zero, which contradicts the assumption of $x^{(1)} \neq x^{(2)}$.
Now we prove that $(2r, 2s)$-sparse strongly observability is a necessary condition. For the sake of contradiction, suppose that the system described by \eqref{sys} is not $(2r, 2s)$-sparse strongly observable, however, reconstructing the state (possibly with delays) is still possible. We construct two system trajectories with different (initial) states that have exactly the same input and output sequences under suitable attack strategies (additive terms). This implies that estimating the correct state is indeed impossible thereby establishing the desired contradiction.\\
By the assumption of the contradiction, the underlying system is not $(2r, 2s)$-sparse strongly observable, so there exist subsets of inputs and outputs denoted by $\Gamma_u $ with $| \Gamma_{u} | = 2r$ and $\Gamma_y $ with $| \Gamma_{y} | = p -2s$, respectively, such that $\subsys{\Gamma_{u}}{ {\Gamma_{y}}}$ is not strongly observable. Corollary \ref{cor:strong_observable1} implies that there exist an initial condition $\Delta x$ and an input sequence $\{ \Delta u(t) \}$ (with its support lying inside $\Gamma_u$) that generates an output sequence $\{ \Delta y(t) \}$ with $\supp{\{ \Delta y(t) \} } \subseteq \overline{\Gamma}_y$. One can rewrite $\Delta u(t)$ and $\Delta y(t)$ as sum of two sparse signals, more precisely:
\begin{align}
\Delta u(t)  &= \Delta u^{(1)}(t) + \Delta u^{(2)}(t), \\
\Delta y(t)  &= \Delta y^{(1)}(t) + \Delta y^{(2)}(t), 
\end{align} where cardinality of $\supp{\{ \Delta u^{(i)}(t) \}}$ and $\supp{\{ \Delta y^{(i)}(t) \}}$ are upper-bounded by $r$ and $s$ for $i \in \{1, 2\}$, respectively. For example, we can rewrite $\overline{\Gamma}_y = \overline{\Gamma}^{(1)}_y \cup \overline{\Gamma}^{(2)}_y$ where $|\overline{\Gamma}^{(i)}_y| \leq s$ for $i \in \{1,2\}$. Then we define 
\begin{align}
\begin{cases}
\Delta y^{(i)}(t)\vert_{\overline{\Gamma}_y^{(i)}} &:= \Delta y(t)\vert_{\overline{\Gamma}^{(i)}} \nonumber \\
\Delta y^{(i)}(t)\vert_{ \Gamma_y^{(i)}} &:= 0
\end{cases}, \quad  \text{for} \quad  i \in \{1,2\}.
\end{align}
Now consider the following two different trajectories of the system
\begin{align}
\begin{cases}
u^{(1)}_S(t) &= \Delta u(t) \\
y^{(1)}_S(t) &= \Delta y(t)
\end{cases}, \quad
\begin{cases}
u^{(2)}_S(t) &= 0 \\
y^{(2)}_S(t) &= 0
\end{cases}
\end{align} with their initial states
\begin{align}
\begin{cases}
x^{(1)}(0) &= \Delta x \\
x^{(2)}(0) &= 0
\end{cases},
\end{align}
 and their corresponding attack strategies,
\begin{align}
\begin{cases}
w^{(1)}(t) &=  \Delta u^{(1)}(t) \\
a^{(1)}(t)  &= -\Delta y^{(1)}(t)
\end{cases}, \quad 
\begin{cases}
w^{(2)}(t) &= -\Delta u^{(2)}(t) \\
a^{(2)}(t)  &= \Delta y^{(2)}(t)
\end{cases}. \label{eq:thm:attack}
\end{align}It is straightforward to verify that $\{y^{(1)}(t)\} = \{ y^{(2)}(t)\}$ and $\{u^{(1)}(t)\} = \{ u^{(2)}(t)\}$, i.e., under the attack model \eqref{model} the controlled inputs and the observed outputs are exactly the same for both trajectories while having different states, therefore the proof is complete.
\end{pf}

\section{Secure Observer Design}
\label{section:secure_observer}

In this section, we seek solutions to Problem \ref{problem:design}. In the first part, we explain the intuition behind the proposed algorithm that estimates the state despite attacks on inputs and outputs. We give formal guarantees that the algorithm reconstructs the state correctly. In the second part, we introduce the observer by leveraging the SMT paradigm followed by two methods that enhance the run time of state estimation.\\
Based on the attack model \eqref{model}, the input to the system is decomposed into two additive terms, the controller-designed input $u(t)$ and the adversarial input $w(t)$. The underlying system \eqref{sys} is linear and therefore we can easily exclude the effect of the controller-designed input from the output by subtracting its effect. Hence, without loss of generality we assume that the true $u(t)$ is zero.\\
The proposed algorithm is based on the following proposition.
\begin{prop}
Suppose the underlying system is $(2r, 2s)$-sparse strongly observable, and the number of attacked inputs and outputs are bounded by $r$ and $s$, respectively. Given any subset of inputs and outputs denoted by ${\Gamma}_u$ and ${\Gamma}_y$ with $|{\Gamma}_u| \leq r$ and $|{\Gamma}_y| \geq p-s$, the first statement below implies the second:
\begin{enumerate}
\item There exist $\bold{\hat{U}} \in \mathbb{R}^{n |T|}$ and $\hat{x} \in \mathbb{R}^{n}$ such that
\begin{align}
\bold{Y}\vert_{{\Gamma}_y}(t) = \mathcal{O}_{{\Gamma}_y} \hat{x} + \mathcal{N}_{{{\Gamma}_u} \to {{\Gamma}_y}} \bold{\hat{U}}.
\label{eq:prop:test}
\end{align}
\item The estimated state $\hat{x}$, is equal to the actual state of the system at time $t-n+1$,  $x(t-n+1)$, where $n$ is the order of the underlying system.
\end{enumerate}
\label{prop:main}
\end{prop}
\begin{rem}
The underlying system is $(2r, 2s)$-sparse strongly observable therefore $\subsys{{\Gamma}_u}{{\Gamma}_y}$ is strongly observable. If \eqref{eq:prop:test} has a solution, then $\hat{x}$ would be the unique solution for $x$ (see section III-B of \cite{partial}).
\end{rem}
\begin{pf}
Let us denote the set of attack-free outputs and under-attack inputs by $\Gamma_y^{\ast}$ and $\Gamma_u^{\ast}$. At most $s$ outputs are under attack, therefore $| {\Gamma}_y \cap \Gamma_y^{\ast} | \geq p-2s $. Note that $\bold{Y}\vert_{{\Gamma}_y \cap \Gamma^{\ast}_y}$ can be written as follows:
\begin{align}
\bold{Y}\vert_{{\Gamma}_y \cap \Gamma_y^{\ast}} =& O_{{\Gamma}_y \cap \Gamma_y^{\ast}} {x(t-n+1)} \nonumber \\ +
&\mathcal{N}_{{\Gamma}_u \to {{\Gamma}_y \cap \Gamma_y^{\ast} }}{ \bold{W}\vert_{{\Gamma}_u}} + \mathcal{N}_{ \Gamma_u^{\ast} \setminus {{\Gamma}_u} \to { {\Gamma}_y \cap \Gamma_y^{\ast}}}{\bold{W}\vert_{\Gamma_u^{\ast} \setminus {{\Gamma}_u}} }.
\end{align} On the other hand, we can rewrite \eqref{eq:prop:test} by taking only outputs in $ {\Gamma}_y \cap \Gamma_y^{\ast}$,
\begin{align}
\bold{Y}\vert_{ {\Gamma}_y \cap \Gamma_y^{\ast}} = O_{ {\Gamma}_y \cap \Gamma_y^{\ast}} {\hat{x}} + \mathcal{N}_{ {\Gamma}_u \to { {\Gamma}_y \cap \Gamma_y^{\ast}}}{\bold{\hat{U}}} + \mathcal{N}_{{\Gamma_u^{\ast} \setminus { {\Gamma}_u}} \to { {\Gamma}_y \cap \Gamma_y^{\ast}}}\bold{0},
\end{align} where $\bold{0}$ is a zero vector with appropriate dimensions. The underlying system is $(2r,2s)$-sparse strongly observable, therefore we conclude that the sub-system $\hat{S} := \subsys{  {\Gamma}_u \cup \Gamma_u^{\ast}  }{ {\Gamma}_y \cap \Gamma_y^{\ast}}$ is strongly observable. One can reinterpret both equations as two (possibly different) valid trajectories of the system $\hat{S}$ that share the same output sequence. Strong observability of $\hat{S}$ implies that $\hat{x} = x(t-n+1) $ which completes the proof.
\end{pf}
The main algorithm in this paper builds upon this proposition. We search for a set of inputs and outputs that satisfies equality \eqref{eq:prop:test}, i.e., we check if there exist $\bold{\hat{U}}$ and $\hat{x}$ that make equality \eqref{eq:prop:test} hold. Based on Proposition \ref{prop:main}, we define a consistency check as follows,
\begin{test}[Consistency Check]
Given subsets of inputs and outputs denoted by ${\Gamma}_u$ and ${\Gamma}_y$, TEST$( {\Gamma}_u, {\Gamma}_y )$ returns true if
\begin{align}
\min_{\bold{\hat{U}}, \hat{x} } \|\bold{Y}\vert_{{\Gamma_y}} - \mathcal{O}_{{\Gamma_y}}\hat{x} - \mathcal{N}_{{\Gamma_u} \to {{\Gamma_y}}} \bold{\hat{U}} \| \leq \epsilon, \label{eq:test}
\end{align} where $\epsilon > 0$ is the solver tolerance, due to numerical errors. However, for the sake of clarity, we focus in this paper on the case when $\epsilon$ is negligible\footnote{Note that the minimum always exists for \eqref{eq:test} as the cost function is a semi-definite quadratic function.}.
\label{test}
\end{test}

Finding the right subset of inputs and ouputs that satisfies this test is a combinatorial problem in nature and requires exhaustive search. It is well-known that secure state estimation under this attack model is in general \emph{NP-hard} \cite{yasser_SMT, fabio}. This test is depicted in Algorithm \ref{alg:t-solver.check}. \\
In the rest of this section, we introduce an architecture for our observer followed by methods to improve its computational performance. For each input (output), we assign a binary variable $\mathbf{b}_i \in \mathbb{B}$ ($\mathbf{c}_i \in \mathbb{B}$) that indicates if the corresponding input (output) is under attack or not, i.e., $\mathbf{b}_i = 1$ ($\mathbf{c}_i = 1$) if the $i^{\text{th}}$ input (output) is under attack. In the rest of this paper, we use the bold letters ($\mathbf{b}$ and $\mathbf{c}$) to denote these Boolean variables and we reserve non-bold type face ($b$ and $c$) as instances of them. Finding the right assignment of these Boolean variables is combinatorial in nature and in order to efficiently decide which set of inputs and outputs satisfies the TEST in \eqref{eq:test}, we design an observer using the lazy SMT paradigm \cite{satisfiability2009}. 

\subsection{Overall Architecture}

\begin{figure}
\centering
\includegraphics[scale=0.6]{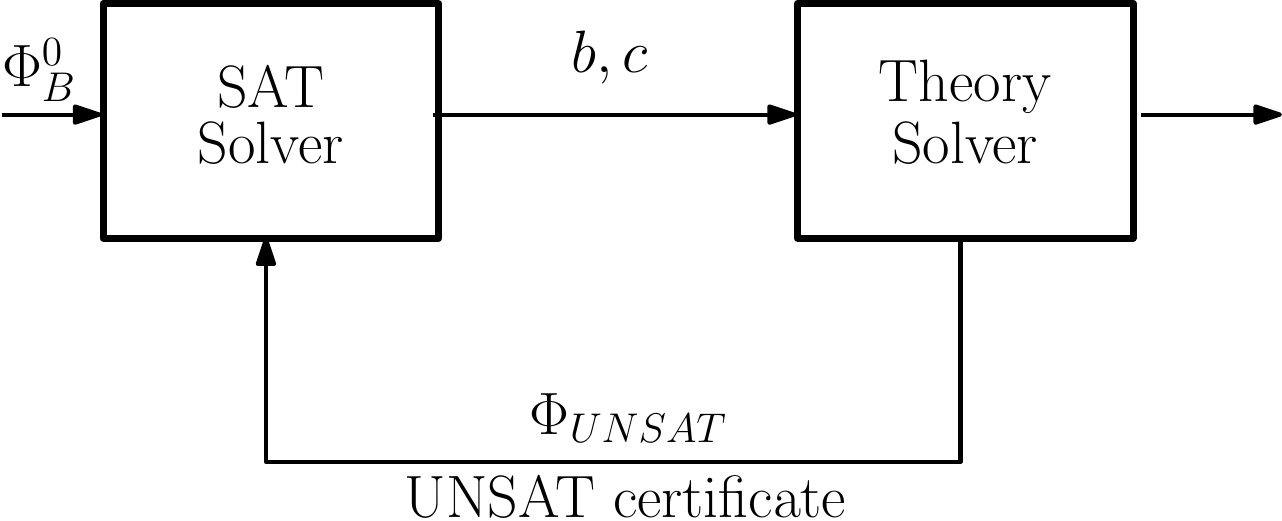}
\caption{The lazy SMT paradigm.}
\label{fig:SMT}
\end{figure}

The observer consists of two blocks that interact with each other, a propositional satisfiability (SAT) solver and a theory solver. The former reasons about the combination of Boolean and pseudo-Boolean constraints and produces a feasible instance of $\mathbf{b} \in \mathbb{B}^{m}$ and $\mathbf{c} \in \mathbb{B}^{p}$, based on its current state. The theory solver checks the consistency of Boolean variables using the consistency test, and when the test fails, it encodes the inconsistency as a pseudo-Boolean constraint and returns it to the SAT solver. The general architecture is depicted in Figure \ref{fig:SMT}.

The initial pseudo-Boolean constraint only bounds the number of attacked inputs and outputs, i.e., 
\begin{align}
\Phi_{B} :=  (\sum\limits_{i = 1}^{m} \mathbf{b}_i \leq r ) \bigwedge (\sum\limits_{j = 1}^{p} \mathbf{c}_i \leq s ).
\end{align} Initially, the SAT solver generates instances of $\mathbf{b}$ and $\mathbf{c}$ that satisfy $\Phi_{B}$. The theory solver checks whether $ {\Gamma_u} := {\supp{b}}$ and ${\Gamma_y} := \overline{\supp{c}}$ satisfies the consistency check. If the test is satisfied, then the algorithm terminates and returns the (delayed) estimate of the state. Otherwise, the theory solver outputs UNSAT and generates a reason for the conflict, a certificate, or a counterexample that is denoted by $\Phi_{\text{cert}}$. This counterexample encodes the inconsistency among the chosen  inputs and outputs. The following always constitutes a naive certificate.
\begin{align}
\Phi_{\text{naive-cert}} := \sum\limits_{i \in \overline{\text{supp}(b)}} \textbf{b}_i + \sum\limits_{j \in \overline{\text{supp}(c)}} \textbf{c}_j  \geq 1. \label{eq:naive}
\end{align}
On the next iteration, the SAT solver updates the constraint by conjoining $\Phi_{\text{cert}}$ to $\Phi_{B}$, and generates another feasible assignment for $\mathbf{b}$ and $\mathbf{c}$. This procedure is repeated until the theory solver returns SAT as illustrated in Algorithm \ref{alg:main}. 

\begin{algorithm}[h]
\caption{Secure state estimator}
\begin{algorithmic}[1]

\Require $A, B, C, D$ (system), $Y$ (output), $r, s$ (bounds)
\State status $\gets$ UNSAT 
\State $\Phi_{\text{cert}} \gets \text{True}$ 
\State $\Phi_{B} \gets (\sum\limits_{i \in \{1, \hdots, m\}} \textbf{b}_i \leq r ) \bigwedge (\sum\limits_{i \in \{1, \hdots, p\}} \textbf{c}_i \leq s )$ 
\While{status $=$ UNSAT}
	\State $\Phi_{B} \gets \Phi_{B} \bigwedge \Phi_{\text{cert}}$
	\State $(b,c) \gets$ SAT-solver($\Phi_B$) 
	\State (status, x) $\gets$ T-solver.check(${\supp{b}},\overline{ \supp{c} }$)
	\State $\Phi_{\text{cert}} \gets $ T-solver.Certificate(${\supp{b}}, \overline{ \supp{c} }$)
\EndWhile
\State {\textbf{return} $(x, b, c)$} 
\end{algorithmic}
\label{alg:main}
\end{algorithm} 

Note that Proposition \ref{prop:main} implies that the SAT solver eventually produces an assignment that satisfies the consistency test and therefore Algorithm \ref{alg:main} always terminates. The size of the certificate plays an important role in the overall execution time of the algorithm \cite{yasser_SMT}. Note that the attack model considered in \cite{yasser_SMT} is restricted to outputs, and the major contribution of our work is to handle both input and output attacks. In the next section, we focus on constructing shorter counterexamples to improve the run time.

\begin{algorithm}
\caption{T-solver.check}
\begin{algorithmic}[1]

\Require {${\Gamma_u}, {\Gamma_y}$ }

\State \textbf{Solve:} $(\hat{x}, \bold{\hat{U}}) = \text{argmin}_{x,\bold{U}} \| \bold{Y}\vert_{{\Gamma_y}} - \mathcal{O}_{{\Gamma_y}}x - \mathcal{N}_{{\Gamma_u} \to {{\Gamma_y}}}\bold{U} \|$
\If{$\|\bold{Y}\vert_{{\Gamma_y}} - \mathcal{O}_{{\Gamma_y}}\hat{x} - \mathcal{N}_{{\Gamma_u} \to {{\Gamma_y}}} \bold{\hat{U}} \| \leq \epsilon $}

	\State status = SAT
\Else

	\State status = UNSAT
\EndIf  
\State {\textbf{return} $(\text{status}, \hat{x} )$}
\end{algorithmic}
\label{alg:t-solver.check}
\end{algorithm}

\subsection{SAT certificate}

In this part, we improve the efficiency of Algorithm \ref{alg:main} by constructing a shorter certificate (counter-example or conflicts). As it was discussed before, the naive certificate only excludes the current assignment of $\mathbf{b}$ and $\mathbf{c}$ from the search space of the SAT solver, however, by exploiting the structure of the underlying system, we show that we can further decrease the size of the certificate and therefore prune the search space more efficiently.\\
One of the main results of this paper is to show that we can always find a smaller conflicting subset of inputs and outputs. We propose two methods for generating shorter certificates. The first method reduces the size of the counterexample by at least $s-1$, we explain this method in Lemma \ref{lemma:sensors} and give a formal proof of the existence of such shorter certificate. In practice, however we observe the reduction in the length of conflicts is much larger than this theoretical bound. The second method is inspired by the \textsc{QuickXplain} algorithm. This method generates counter-examples that are irreducible, meaning that we cannot reduce the size of the counter-example by removing some of it's entries. We also note that by generating multiple certificates at each iteration we can further enhance the execution time. At the end of this section Lemma \ref{lemma:hueristics} states that for a generic LTI system the size of the certificate cannot be smaller than $m+1$.\\
Let us assume that the SAT solver hypothesized {$\Gamma_u^{\text{SAT}} := {\supp{b}}$} and \mbox{$\Gamma_y^{\text{SAT}} := \overline{\supp{c}}$} as the set of compromised inputs and safe outputs, respectively. The main intuition behind both methods is to look for $\Gamma_u^{\text{cert}} \supseteq \Gamma_u^{\text{SAT}}$ and \mbox{$\Gamma_y^{\text{cert}} \subseteq \Gamma_y^{\text{SAT}} $} that would not satisfy the consistency test. Note that the certificate consists of inputs in $\overline{\Gamma}_u^{\text{cert}}$ and outputs in $\Gamma_y^{\text{cert}}$.

\subsection{Method I based on heuristics}

Method I reduces the size of the certificate by increasing the size of (supposedly under attack) inputs ($\Gamma_u^{\text{cert}}$) followed by decreasing the size of (supposedly safe) outputs ($\Gamma_u^{\text{cert}}$). The summary of the above procedure of shortening certificates is illustrated in Algorithm \ref{alg:t-solver.certificate1}. We begin by adding inputs to $\Gamma_u^{\text{SAT}}$ while making sure TEST still returns false and the number of inputs is bounded by $2r$. Let us denote this new set of inputs by $\Gamma_u^{\text{cert}}$.\\
At the second step, we shrink the set of conflicting outputs in order to further shorten the size of the counterexample. Let us denote a subset of $\Gamma_y^{\text{SAT}}$ of size $p-2s$ by $\Gamma_y^\text{temp}$. The following lemma shows we can reduce the size of conflicting outputs at least by $s-1$.

\begin{algorithm}
\caption{T-solver.Certificate 1}

\begin{algorithmic}[1]

\Require {$\Gamma_u^{\text{SAT}}, \Gamma_y^{\text{SAT}} $}

\Statex \textbf{step 1:} Conduct a linear search in the input set 
\State Sort $\overline{\Gamma}_u^{\text{SAT}}$ \;
\State $\text{status} \gets \text{UNSAT}, j \gets \emptyset, \Gamma_u^{\text{cert}} \leftarrow  \Gamma_u^{\text{SAT}}$ \;
\While{ status == UNSAT {\textbf{and}} $|\Gamma_u^{\text{cert}}| < 2r$ }
	\State $\Gamma_u^{\text{cert}} \leftarrow \Gamma_u^{\text{cert}} \cup \{ j \}$
	\State pick another input $j \in \overline{\Gamma}_u^{\text{SAT}}$
	\State $(\text{status}, x) \leftarrow$ T-Solver.check$(\Gamma_u^{\text{cert}} \cup \{ j \}, \Gamma_y^{\text{SAT}})$
	
\EndWhile
\Statex \textbf{step 2:} Conduct a linear search in the output set
\State Sort $\Gamma_y^{\text{SAT}}$ 
\State Pick a subset of size $p-2s$: $\Gamma_y^{\text{temp}} \subseteq \Gamma_y^{\text{SAT}}$
\State $\text{status} \gets \text{SAT}, i \gets \emptyset$ 
\While{ status == SAT }
	\State $\Gamma_y^{\text{cert}}  \leftarrow \Gamma_y^{\text{temp}}  \cup \{i\}$
	\State $(\text{status}, x) \gets$ T-Solver.check$(\Gamma_u^{\text{cert}}, \Gamma_y^{\text{cert}})$
	\State Pick another output  $i \in  \Gamma_y^{\text{SAT}} \setminus \Gamma_y^{\text{temp}}  $
\EndWhile
\State $\Phi_{\text{cert}}^1 \gets \sum\limits_{j \in \overline{\Gamma}_u^{\text{cert}}} \mathbf{b}_j + \sum\limits_{i \in {\Gamma}_y^{\text{cert}} } \mathbf{c}_i \geq 1$
\State {\textbf{return} $\Phi_{\text{cert}}^1$ }

\end{algorithmic}
\label{alg:t-solver.certificate1}
\end{algorithm}




\begin{lem}
Assume that the system $S$ is $(2r, 2s)$-sparse strongly observable, and the number of attacked inputs and outputs are bounded by $r$ and $s$, respectively. Pick any subset of inputs and outputs denoted by $\Gamma_u^{\text{cert}}$ and $\Gamma_y^{\text{SAT}}$ with $|\Gamma_u^\text{cert}| \leq 2r$ and $|\Gamma_y^{\text{SAT}}| \geq p-s$, that do not satisfy the consistency check \eqref{eq:test}. Given any subset of at most $p-2s$ outputs denoted by $\Gamma_y^{\text{temp}} \subseteq \Gamma_y^{\text{SAT}}$, one of the following is true:
\begin{enumerate}
\item TEST($\Gamma_u^\text{cert}, \Gamma_y^\text{temp}$) returns false,
\item There exists an output $i \in \Gamma_y^{\text{SAT}} \setminus \Gamma_y^\text{temp}$ such that TEST($\Gamma_u^\text{cert}, \Gamma_y^\text{temp} \cup \{i\}$) returns false.
\end{enumerate} \label{lemma:sensors}
\end{lem}
\begin{pf}
Please refer to Appendix.
\end{pf}

We denote this smaller set of conflicting outputs $\Gamma_y^\text{temp}$(if TEST($\Gamma_u^\text{cert}, \Gamma_y^{\text{temp}} $) returns false, otherwise $\Gamma_y^\text{temp} \cup \{ i \}$) by $\Gamma_y^\text{cert}$. {Lemma \ref{lemma:sensors} gives formal guarantees of the existence of shorter certificates which hold no matter how the subsets of inputs and outputs ($\Gamma_u^\text{temp}$ and $\Gamma_y^\text{temp}$) are chosen. This lemma shows that Method I reduces the size of the certificate by at least $s-1$.\\
In practice, we choose these subsets based on heuristics that have for objective a decrease in the overall running time. We assign slack variables to inputs and outputs similarly to \cite{yasser_SMT} and \cite{showkatbakhsh3}, and sort them based on the structure of the system}. Recall that Algorithm \ref{alg:t-solver.certificate1} shortens the certificate by reducing the number of inputs followed by the reduction in the number of outputs, i.e., we \emph{simultaneously} reducing both inputs and outputs in the certificate. We observe that by generating two counterexamples, we can prune the search space of the SAT solver more efficiently. Similarly to Algorithm \ref{alg:t-solver.certificate2}, we can find two counterexamples by reducing the number of inputs following a reduction in the number of outputs and vice-verse.\\
\noindent\textbf{Sorting $\overline{\Gamma}_u^{\text{SAT}}$ and $ {\Gamma}_y^{\text{SAT}}$:}\\
Assuming TEST($\Gamma_u^{\text{SAT}}, \Gamma_y^{\text{SAT}}$) returns false, we assign slack variables to inputs in $\overline{\Gamma}_u^{\text{SAT}}$ and outputs in $\Gamma_y^{\text{SAT}}$, denoted by $\text{slack}_{u}(j)$ and $\text{slack}_{y}(i)$, respectively. Let us denote a solution to the optimization \eqref{eq:test} inside TEST($\Gamma_u^{\text{SAT}}, \Gamma_y^{\text{SAT}}$) by $\hat{x}$ and $\hat{\bold{U}}$.\\
We define $\text{slack}_{u}(j)$ for $j \in  \overline{\Gamma}_u^{\text{SAT}}$ as the norm of the projection of $\bold{Y}\vert_{{\Gamma_y^{\text{SAT}}}} - \mathcal{O}_{{\Gamma_y^{\text{SAT}}}}\hat{x} - \mathcal{N}_{{\Gamma_u^{\text{SAT}}} \to {{\Gamma_y^{\text{SAT}}}}} \bold{\hat{U}}$ onto the column space of $\mathcal{N}_{j \to {\Gamma}_y^{\text{SAT}}}$,
\begin{align}
\text{slack}_{u}(j) :=& \\
\| \mathcal{N}_{j \to {\Gamma}_y^{\text{SAT}}} & \mathcal{N}_{j \to {\Gamma}_y^{\text{SAT}}}^{\dagger}  \left( \bold{Y}\vert_{{\Gamma_y^{\text{SAT}}}} - \mathcal{O}_{{\Gamma_y^{\text{SAT}}}}\hat{x} - \mathcal{N}_{{\Gamma_u^{\text{SAT}}} \to {{\Gamma_y^{\text{SAT}}}}} \bold{\hat{U}} \right) \|. \nonumber
\end{align} This slack variable measures how much of the residual can be justified by considering $j$ in addition to ${\Gamma}_u^{\text{SAT}}$. Note that we want to append inputs to ${\Gamma}_u^{\text{SAT}}$ while having a false TEST. We first normalize these slack variables by the norm of the corresponding invertibility matrix, and $\overline{\Gamma}_u^{\text{SAT}} $ is obtained by sorting slack variables in \emph{ascending} order.\\
We define $\text{slack}_{y}(i)$ as the residual of each output:
\begin{align}
\text{slack}_{y}(i) := \| \bold{Y}\vert_{i} -  \mathcal{O}_{i} \hat{x} - \mathcal{N}_{ \Gamma_u^{\text{SAT}} \to \{ i \} } {\bold{U} } \|, \quad i \in \Gamma_y^{\text{SAT}}. \label{eq:heuristics:sensors}
\end{align} Note that,
\begin{align}
\sum\limits_{i \in \Gamma_u^{\text{SAT}}} \text{slack}_{y}(i) = \min_{\bold{\hat{U}}, \hat{x} } \| \bold{Y}\vert_{{\Gamma_y^{\text{SAT}}}} - \mathcal{O}_{{\Gamma_y^{\text{SAT}}}}\hat{x} - \mathcal{N}_{{\Gamma_u^{\text{SAT}}} \to {{\Gamma_y^{\text{SAT}}}}} \bold{\hat{U}} \|.
\end{align}
We first normalize each slack variable by the norm of the corresponding observabality matrix. Recall that we aim to find a smaller subset of $\Gamma_u^{\text{SAT}}$ while ensuring TEST returns false. We pick the output with the highest slack variable as the first element of $\Gamma_u^{\text{SAT}}$. We sort the rest based on the dimension of the kernel of each observability matrix, following the intuition provided in \cite{yasser_SMT}.
\subsection{Method II based on QuickXplain}

The second method (Algorithm \ref{alg:t-solver.certificate2}) is inspired by \textsc{QuickXplain} and generates a counter-example by pruning the naive-certificate \eqref{eq:naive} to make it irreducible. We formally define this property as follows,
\begin{defn}[Irreducible certificate]
A certificate consisting of inputs $\overline{\Gamma}_u$ and outputs $\Gamma_y$ is irreducible, if no other subset of it can generate a conflict, i.e., for all subsets denoted by $ \overline{\Gamma}'_u \subseteq \overline{\Gamma}_u$ and $\Gamma'_y \subseteq \Gamma_y$ the following are equivalent:
\begin{enumerate}
\item $ \overline{\Gamma}'_u$ and $\Gamma'_y$ generate a conflict. 
\item $ \overline{\Gamma}'_u = \overline{\Gamma}_u$ and $\Gamma'_y = \Gamma_y$.
\end{enumerate}
\end{defn}
One cannot prune irreducible certificates and each element is necessary for the set to remain a counter-example. Let $\Delta^{\text{SAT}} $ be the elements (consisting of inputs $\overline{\Gamma}_u^{\text{SAT}}$ and outputs $\Gamma_y^{\text{SAT}}$) of the naive certificate. For ease of exposition we slightly abuse notation to denote $\mbox{TEST}( \Gamma_u^{\text{SAT}}, \Gamma_y^{\text{SAT}} )$ by $\mbox{TEST}(\Delta^{\text{SAT}})$. We denote the output of this algorithm by $\Delta_{\text{cert}}$ which consists of inputs $\overline{\Gamma}_u^{\text{cert}}$ and outputs $\Gamma_y^{\text{cert}}$.\\
This method consists of an exploration phase in which it finds an element (input or output) that belongs to an irreducible certificate. Let us denote an enumeration of $\Delta^{\text{SAT}}$ by $e_1, \cdots, e_k$, and the internal state by $\Delta_{\text{temp}} \leftarrow \emptyset$. This method begins by adding step-by-step elements of $\Delta^{\text{SAT}}$ to $\Delta_{\text{temp}}$. The first element ($e_i \in \Delta^{\text{SAT}}$) that fails TEST$(\Delta_{\text{temp}})$ is part of an irreducible certificate, and therefore is added to $\Delta_{\text{cert}}$.\\
In order to find further elements of this certificate, we keep $e_i$ in the background and the first element that fails the consistency check is added to $\Delta_{\text{cert}}$. This repeated process can be implemented efficiently by using the divide and conquer paradigm as depicted in Algorithm \ref{alg:t-solver.quickxplain}. When an element $e_i$ of $\Delta^{\text{SAT}}$ is detected we divide the the remaining elements into two disjoint subsets $\Delta^1 := \{ e_1, \cdots, e_j \}$ and $\Delta^2 := \{ e_{j+1}, \cdots, e_{i-1} \}$. We can now recursively apply the algorithm to find a conflict $\Delta_{\text{cert}}^2$ among $\Delta^2$ by keeping the set $\Delta^1$ in the background and a conflict $\Delta_{\text{cert}}^1$ among $\Delta^1$ by keeping the set $\Delta_{\text{cert}}^2$ in the background. This method of finding an irreducible subset is depicted in Algorithm \ref{alg:t-solver.quickxplain} \\
Note that the resulting counter-example depends on the initial enumeration of elements in $\Delta^{\text{SAT}}$. If the all the inputs (outputs) are ahead of outputs (inputs), then the resulting counter-example mostly consists of inputs (outputs). In order to have the maximal reduction in the search space of the SAT solver at each iteration, we produce three certificate using this method, putting inputs first, outputs first and mixing both inputs and outputs.\\
In the last part of this section, we look at the certificate size for a generic LTI system. We observe that the certificate size cannot be smaller that the number of inputs which is stated formally in the following lemma.

\begin{algorithm}
\caption{T-solver.QuickXplain}
\begin{algorithmic}[1]

\Require {$\Delta_{\text{cert}}^0, \Delta^0 $ }

\If{T-solver.check($\Delta_{\text{cert}}^0$) = UNSAT or $\Delta^0 == \emptyset$}  
	\State \textbf{return} $\emptyset$ 
\EndIf
\Statex Let $e_1, \cdots, e_k$ be an enumeration of $\Delta^0$ 
\State $i \leftarrow 0$, $\Delta_{\text{temp}} \leftarrow \Delta_{\text{cert}}^0,$
\While{ T-solver.check($\Delta_{\text{temp}}$) = SAT and $i \leq k$}

	\State $i \gets i+1$ 
	\State $\Delta_{\text{temp}} \gets \Delta_{\text{temp}} \cup e_i$ 
	\State $\Delta_{\text{cert}}^i \leftarrow \Delta_{\text{temp}}$ 

\EndWhile

\State $\Delta_{\text{cert}} \leftarrow e_{i}$, $j \leftarrow \lfloor \frac{i}{2} \rfloor$
\State $\Delta^1 \leftarrow \{ e_1, \cdots, e_j \}$ 
\State $\Delta^2 \leftarrow \{ e_{ j +1 }, \cdots, e_{i-1} \}$ 

\State $\Delta_{\text{cert}} \leftarrow$ \mbox{$\Delta_{\text{cert}} \cup$ T-solver.QuickXplain($\Delta_{\text{cert}}^{j} \cup \Delta_{\text{cert}}, \Delta^2$)}
\State $\Delta_{\text{cert}} \leftarrow \Delta_{\text{cert}} \cup$ T-solver.QuickXplain($\Delta_{\text{cert}}^{0}\cup\Delta_{\text{cert}}, \Delta^1$)
\State {\textbf{return} $\Delta_{\text{cert}}$ }
\end{algorithmic}
\label{alg:t-solver.quickxplain}
\end{algorithm} 

\begin{algorithm}
\caption{T-solver.Certificate 2}
\begin{algorithmic}[1]

\Require {$\Gamma_u^{\text{SAT}}, \Gamma_y^{\text{SAT}} $ }

\State $\Delta_{\text{cert}} \leftarrow$ T-solver.QuickXplain($\emptyset ,\overline{\Gamma}_u^{\text{SAT}} \cup \Gamma_y^{\text{SAT}}$)
\State Divide $\Delta_{\text{cert}}$ to inputs $\overline{\Gamma}_u^{\text{cert}}$ and outputs $\Gamma_y^{\text{cert}}$
\State $\Phi_{\text{cert}}^2 \leftarrow \sum\limits_{j \in \overline{\Gamma}_u^{\text{cert}}} \mathbf{b}_j + \sum\limits_{i \in \Gamma_y^{\text{cert}} } \mathbf{c}_i \geq 1$
\State {\textbf{return} $\Phi_{\text{cert}}^2$ }
\end{algorithmic}
\label{alg:t-solver.certificate2}
\end{algorithm} 
\begin{lem}
For a generic LTI system the size of the certificate is always lower bounded by $m+1$, where $m$ is the number of inputs. \label{lemma:hueristics}
\end{lem}
\begin{pf}
Please refer to Appendix.
\end{pf}
\section{Simulation Results}

\label{section:simulation_results}

We implemented our SMT-based estimator in \textsc{Matlab} while interfacing with the SAT solver SAT4J \cite{sat4j} and assessed its performance in two case studies, randomly generated LTI systems and a chemical plant. We report the overall running time by using the two proposed methods, Algorithm \ref{alg:t-solver.certificate1} and Algorithm \ref{alg:t-solver.certificate2}.
\subsection{Random Systems}

\begin{figure}
\centering
\includegraphics[scale=0.4]{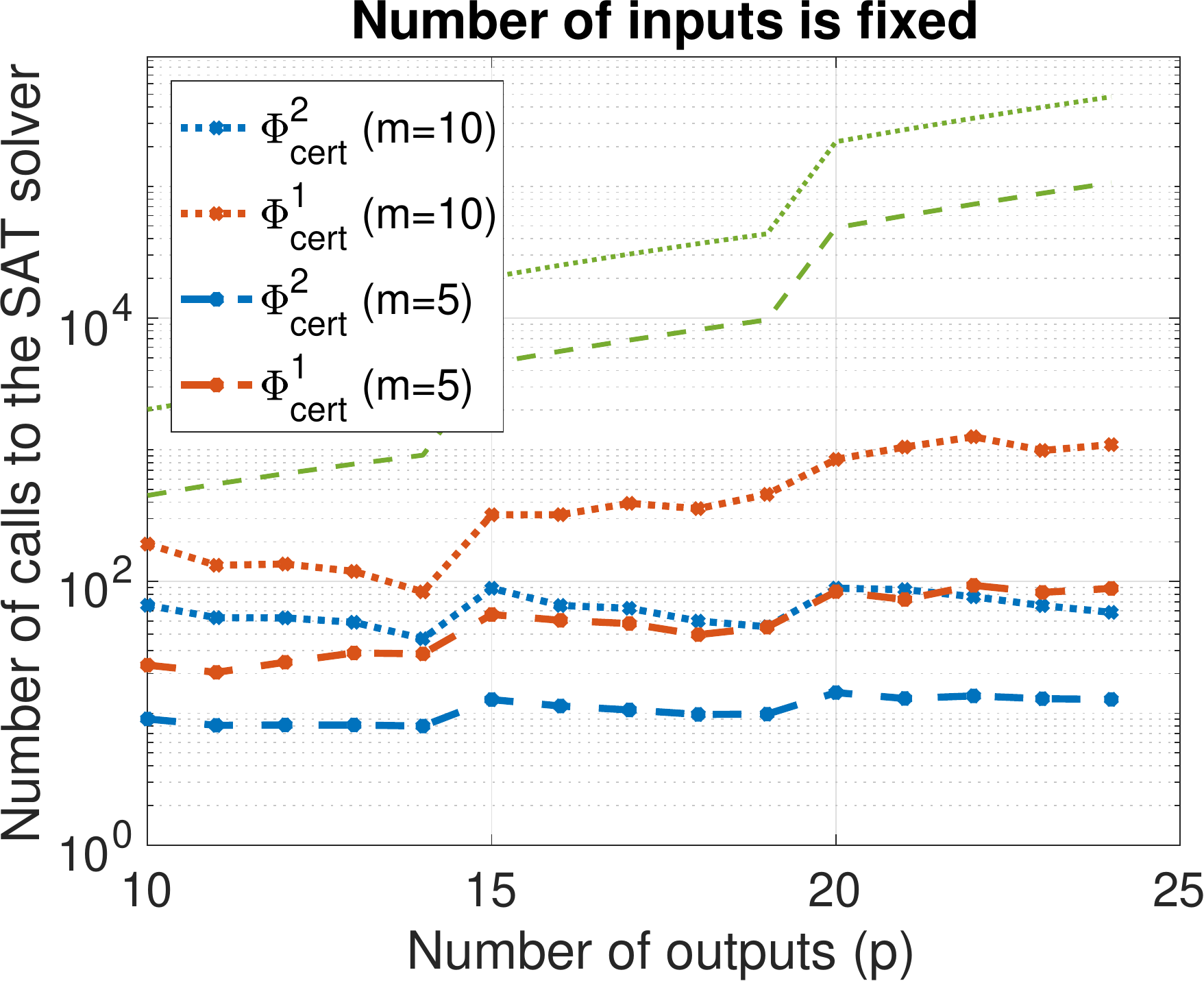}
\caption{Number of calls to the SAT solver in Algorithm \ref{alg:main} using $\Phi_{\text{cert}}^1$, $\Phi_{\text{cert}}^2$ versus the number of outputs ($p$) for a fixed number of inputs. Green dotted and green dashed lines are upper-bounds for the number of the SAT solver calls when using the naive certificate for $m=5$ and $m=10$, respectively. } 
\label{fig:1_}
\end{figure}
\begin{figure}
\centering
\includegraphics[scale=0.4]{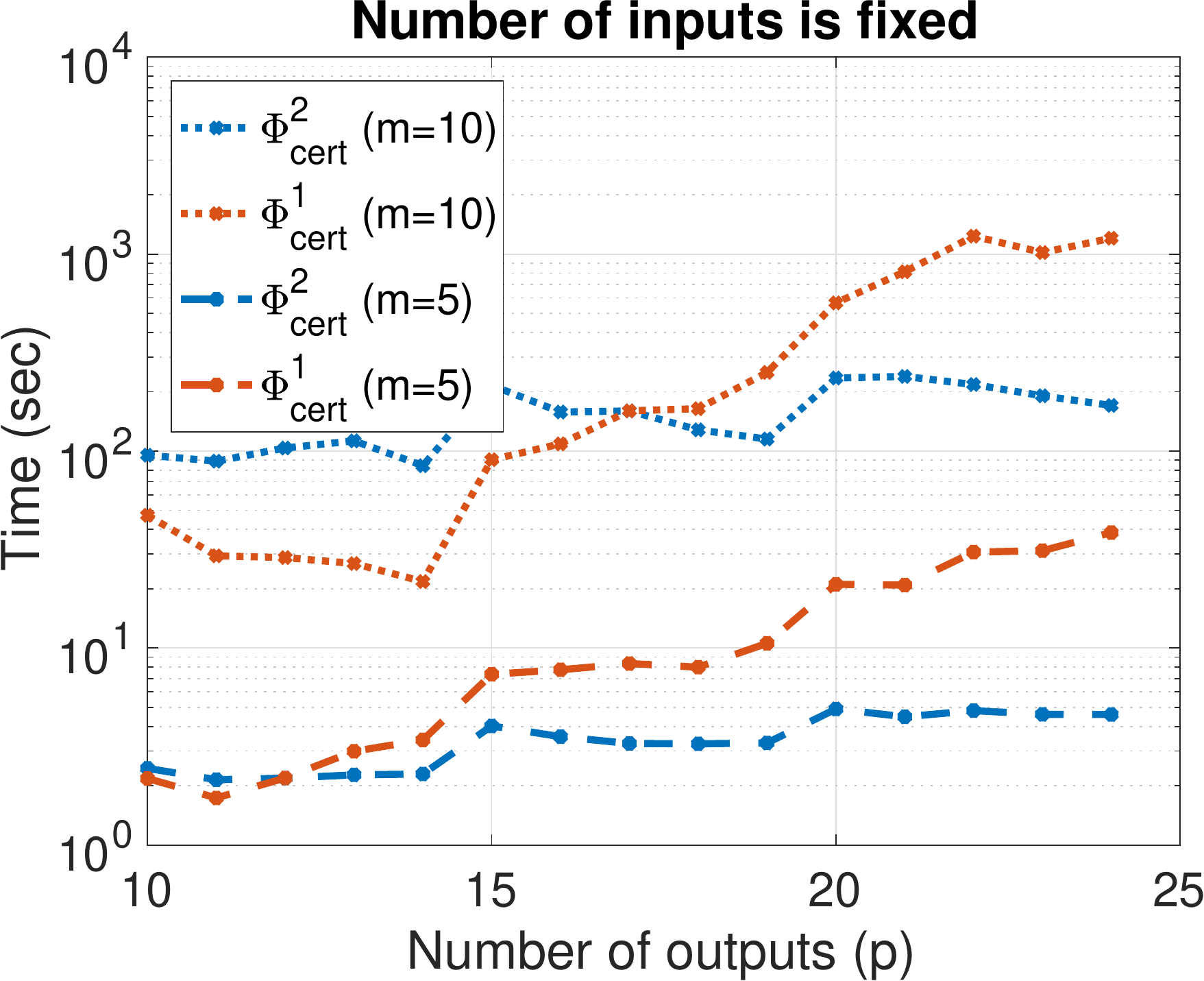}
\caption{Execution time of Algorithm \ref{alg:main} using $\Phi_{\text{cert}}^1$, $\Phi_{\text{cert}}^2$ versus the number of outputs ($p$) and inputs ($m$).  } 
\label{fig:2_}
\end{figure}

We randomly generate systems with a fixed state dimension ($n = 40$) and increase the number of inputs and outputs. Each system is generated by drawing entries of $(A, B, C, D)$ according to uniform distribution, when necessary we scale $A$ to ensure that the spectral radius is close to one. In each experiment, twenty percent of inputs and outputs are under adversarial attacks, and we generate the support set for the adversarial signals uniformly at random. Attack signals and the initial states are drawn according to independent and normally distributed random variables with zero mean and unit variance.
All the systems under experiment satisfy a suitable sparse strong observability condition as described in Section \ref{section:conditions}. \\
Figures \ref{fig:1_} and \ref{fig:2_} report the results of the simulations, each point represents the average of $20$ experiments. All the experiments run on an Intel Core i5 2.7GHz processor with 16GB of RAM. We verify the run-time improvement resulting from using the shorter certificates, $\Phi_{\text{cert}}^1$ and $\Phi_{\text{cert}}^2$, compared to the theoretical upper-bound of the brute-force approach in Figure \ref{fig:1_}. For instance, consider the scenario with $p = 24$ and $m = 10$ in Figures \ref{fig:1_} and \ref{fig:2_}. In the brute-force approach, we require to check all $24 \choose 4$ $\times$ $10 \choose 2$ $\approx 4.8*10^{5} $ different combinations of inputs and outputs, however, by exploiting either $\Phi_{\text{cert}}^1$ or $\Phi_{\text{cert}}^2$ we observe a substantial improvement. 
We observe that although $\Phi_{\text{cert}}^2$ gives a worse run time for systems with smaller number of outputs, it scales better compared to $\Phi_{\text{cert}}^1$ when the number of inputs and outputs grow.

\subsection{Chemical Plant}
In this part, we use the proposed observer to detect attacks on inputs and outputs of a simplified version of the Tennessee Eastman control challenge problem \cite{plant}. Ricker \cite{ricker93} derived a continuous time LTI model of the plant interaction in its steady state. 
This system consists of $4$ control inputs and $10$ measured outputs and the linearized model has $8$ state variables. The structure of the continuous-time dynamics is reported below.
\begin{align*}
\frac{dx}{dt} &= {\tiny \begin{bmatrix}
\ast & \ast & \ast & \ast & \ast & \ast & \ast & 0 \\
\ast & \ast & \ast & \ast & \ast & 0 & \ast & 0 \\
\ast & \ast & \ast & \ast & \ast & 0 & \ast & 0 \\
\ast & \ast & \ast & \ast & 0 & 0 & 0 & \ast \\
0 & 0 & 0 & 0 & \ast & 0 & 0 & 0 \\
0 & 0 & 0 & 0 & 0 & \ast & 0 & 0 \\
0 & 0 & 0 & 0 & 0 & 0 & \ast & 0 \\
0 & 0 & 0 & \ast & 0 & 0 & 0 & \ast \\
\end{bmatrix} }x + {\tiny \begin{bmatrix}
0 & 0 & 0 & 0 \\
0 & 0 & 0 & 0 \\
0 & 0 & 0 & 0 \\
0 & 0 & 0 & 0 \\
\ast & 0 & 0 & 0 \\
0 & \ast & 0 & 0 \\
0 & 0 & \ast & 0 \\
0 & 0 & 0 & \ast \\
\end{bmatrix} } u, \\
y &= {\tiny \begin{bmatrix}
0 & 0 & 0 & 0 & \ast & 0 & 0 & 0 \\
0 & 0 & 0 & 0 & 0 & \ast & 0 & 0 \\
\ast & \ast & \ast & \ast & 0 & 0 & \ast & 0 \\ 
\ast & \ast & \ast & \ast & 0 & 0 & 0 & \ast \\
\ast & \ast & \ast & \ast & 0 & 0 & 0 & 0 \\
0 & 0 & 0 & \ast & 0 & 0 & 0 & 0 \\
\ast & \ast & \ast & 0 & 0 & 0 & 0 & 0 \\
\ast & \ast & \ast & 0 & 0 & 0 & 0 & 0 \\
\ast & \ast & \ast & 0 & 0 & 0 & 0 & 0 \\
\ast & \ast & \ast & 0 & 0 & 0 & \ast & \ast \\
\end{bmatrix} }x,
\end{align*} where $\ast$ represents a non-zero entry\footnote{For the exact dynamics of the LTI model, see \cite{ricker93}}, and $x \in \mathbb{R}^8$, $u \in \mathbb{R}^4$ and $y \in \mathbb{R}^{10}$ are state, input and output variables, respectively.
 The only known limitation of this LTI model is the system should operate close to its steady-state. We obtain a discrete-time model by discretizing the continuous-time model assuming a zero-order hold for the input $u$, with a time-step of $5s$. The attacker can read all the inputs and outputs and manipulate one control input and two measured outputs. The linearized system is $(2,4)$-sparse strongly observable, therefore our observer can correctly reconstruct the state under this attack model.\\
We randomly generate attack signals and the initial state according to independent and normally distributed random variables. The support set of attacks are drawn uniformly at random, and in each experiment one input and two outputs are under adversarial attacks.\\
The proposed observer in this paper can correctly reconstruct the (delayed) state after $8$ samples, and the {average} performance {of $20$ experiments}, by using $\Phi_{\text{cert}}^1$ and $\Phi_{\text{cert}}^2$ is reported in Table \ref{table1:sim}. The overall execution time is the run time of the observer after receiving all the required samples from the plant, and it does not take the sampling time of the plant into account. We observe that the execution time of the observer to reconstruct the state and to detect attacks is much smaller compared to the sampling time of the plant.

\begin{table}
\centering
\caption{Average performance of the proposed observer}
\begin{tabular}{c|c|c}
& Overall execution time & \begin{tabular}{c} Number of calls to \\ the SAT solver \end{tabular} \\ \hline
$\Phi_{\text{cert}}^1$ & $0.22 s$ & $20.05$ \\  \hline
$\Phi_{\text{cert}}^2$ & $0.21 s$ & $7.95$
\end{tabular}
\label{table1:sim}
\end{table}

\section{Conclusion}
\label{section:conc}
In this paper, we considered the problem of secure state estimation when inputs and/or outputs are under adversarial attacks. In this set-up, there is no restriction on how the adversary manipulates inputs and outputs. By introducing the notion of sparse strong observability, we derived necessary and sufficient conditions under which state estimation is possible given bounds on the number of attacked outputs and inputs. Furthermore, we demonstrated the scalability and effectiveness of the proposed estimator with numerical simulations.

\bibliographystyle{plain}        
\bibliography{references}

\begin{thebibliography}{10}

\bibitem{par1_1}
Saurabh Amin, Galina~A Schwartz, and Amir Hussain.
\newblock In quest of benchmarking security risks to cyber-physical systems.
\newblock {\em IEEE Network}, 27(1):19--24, 2013.

\bibitem{gupta2}
Cheng-Zong Bai, Vijay Gupta, and Fabio Pasqualetti.
\newblock On kalman filtering with compromised sensors: Attack stealthiness and
  performance bounds.
\newblock {\em IEEE Trans. Autom. Control}, 62(12):6641--6648, 2017.

\bibitem{gupta1}
Cheng-Zong Bai, Fabio Pasqualetti, and Vijay Gupta.
\newblock Data-injection attacks in stochastic control systems: Detectability
  and performance tradeoffs.
\newblock {\em Automatica}, 82:251--260, 2017.

\bibitem{satisfiability2009}
Clark~W Barrett, Roberto Sebastiani, Sanjit~A Seshia, and Cesare Tinelli.
\newblock Satisfiability modulo theories.
\newblock {\em Handbook of satisfiability}, 185:825--885, 2009.

\bibitem{faulttolerant}
Mogens Blanke, Michel Kinnaert, Jan Lunze, Marcel Staroswiecki, and
  J~Schr{\"o}der.
\newblock {\em Diagnosis and fault-tolerant control}, volume 691.
\newblock Springer, 2006.

\bibitem{cardenas}
Alvaro~A C{\'a}rdenas, Saurabh Amin, and Shankar Sastry.
\newblock Research challenges for the security of control systems.
\newblock In {\em HotSec}, 2008.

\bibitem{chong}
Michelle~S Chong, Masashi Wakaiki, and Joao~P Hespanha.
\newblock Observability of linear systems under adversarial attacks.
\newblock In {\em American Control Conference (ACC)}, pages 2439--2444, 2015.

\bibitem{paulo_persis1}
Claudio De~Persis and Pietro Tesi.
\newblock Input-to-state stabilizing control under denial-of-service.
\newblock {\em IEEE Transactions on Automatic Control}, 60(11):2930--2944,
  2015.

\bibitem{plant}
James~J Downs and Ernest~F Vogel.
\newblock A plant-wide industrial process control problem.
\newblock {\em Computers \& chemical engineering}, 17(3):245--255, 1993.

\bibitem{robust2}
S.~Farahmand, G.~B. Giannakis, and D.~Angelosante.
\newblock Doubly robust smoothing of dynamical processes via outlier sparsity
  constraints.
\newblock {\em IEEE Transactions on Signal Processing}, 59(10):4529--4543, Oct
  2011.

\bibitem{fawzi}
Hamza Fawzi, Paulo Tabuada, and Suhas Diggavi.
\newblock Secure estimation and control for cyber-physical systems under
  adversarial attacks.
\newblock {\em IEEE Transactions on Automatic Control}, 59(6):1454--1467, 2014.

\bibitem{ACMsurvey18}
Jairo Giraldo, David Urbina, Alvaro Cardenas, Junia Valente, Mustafa Faisal,
  Justin Ruths, Nils~Ole Tippenhauer, Henrik Sandberg, and Richard Candell.
\newblock A survey of physics-based attack detection in cyber-physical systems.
\newblock {\em ACM Computing Surveys (CSUR)}, 51(4):76, 2018.

\bibitem{jeep}
Andy Greenberg.
\newblock Hackers remotely kill a jeep on the highway, with me in it.
\newblock {\em [online]
  http://www.wired.com/2015/07/hackers-remotely-kill-jeep-highway}, 2015.

\bibitem{basar1}
Abhishek Gupta, C{\'e}dric Langbort, and Tamer Basar.
\newblock Optimal control in the presence of an intelligent jammer with limited
  actions.
\newblock In {\em 49th IEEE Conference on Decision and Control (CDC)}, pages
  1096--1101, 2010.

\bibitem{necmiye}
Farshad Harirchi and Necmiye Ozay.
\newblock Guaranteed model-based fault detection in cyber-physical systems: A
  model invalidation approach.
\newblock {\em arXiv preprint arXiv:1609.05921}, 2016.

\bibitem{strong_obs}
M.L.J. Hautus.
\newblock Strong detectability and observers.
\newblock {\em Linear Algebra and its Applications}, 50(Supplement C):353 --
  368, 1983.

\bibitem{jones}
Harold~Lee Jones.
\newblock {\em Failure detection in linear systems.}
\newblock PhD thesis, Massachusetts Institute of Technology, 1973.

\bibitem{quickxplain}
Ulrich Junker.
\newblock Quickxplain: Conflict detection for arbitrary constraint propagation
  algorithms.
\newblock In {\em IJCAI’01 Workshop on Modelling and Solving problems with
  constraints}, 2001.

\bibitem{nissan}
Leo Kelion.
\newblock Nissan leaf electric cars hack vulnerability disclosed.
\newblock {\em [online] http://www.bbc.com/news/technology-35642749}, 2016.

\bibitem{stuxnet}
Ralph Langner.
\newblock Stuxnet: Dissecting a cyberwarfare weapon.
\newblock {\em IEEE Security \& Privacy}, 9(3):49--51, 2011.

\bibitem{sat4j}
Daniel Le~Berre and Anne Parrain.
\newblock The sat4j library, release 2.2, system description.
\newblock {\em Journal on Satisfiability, Boolean Modeling and Computation},
  7:59--64, 2010.

\bibitem{robust1}
J.~Mattingley and S.~Boyd.
\newblock Real-time convex optimization in signal processing.
\newblock {\em IEEE Signal Processing Magazine}, 27(3):50--61, May 2010.

\bibitem{shaunak}
Shaunak Mishra, Yasser Shoukry, Nikhil Karamchandani, Suhas Diggavi, and Paulo
  Tabuada.
\newblock Secure state estimation: Optimal guarantees against sensor attacks in
  the presence of noise.
\newblock {\em IEEE Transactions on Control of Network Systems}, 4(1):49--59,
  2017.

\bibitem{mo14}
Yilin Mo, Rohan Chabukswar, and Bruno Sinopoli.
\newblock Detecting integrity attacks on scada systems.
\newblock {\em IEEE Transactions on Control Systems Technology},
  22(4):1396--1407, 2014.

\bibitem{yilin1}
Yilin Mo, Emanuele Garone, Alessandro Casavola, and Bruno Sinopoli.
\newblock False data injection attacks against state estimation in wireless
  sensor networks.
\newblock In {\em 49th IEEE Conference on Decision and Control (CDC)}, pages
  5967--5972, 2010.

\bibitem{par1_2}
Yilin Mo, Tiffany Hyun-Jin Kim, Kenneth Brancik, Dona Dickinson, Heejo Lee,
  Adrian Perrig, and Bruno Sinopoli.
\newblock Cyber--physical security of a smart grid infrastructure.
\newblock {\em Proceedings of the IEEE}, 100(1):195--209, 2012.

\bibitem{yilin3}
Yilin Mo and Bruno Sinopoli.
\newblock Secure control against replay attacks.
\newblock In {\em 47th Annual Allerton Conference on Communication, Control,
  and Computing}, pages 911--918. IEEE, 2009.

\bibitem{mo16}
Yilin Mo and Bruno Sinopoli.
\newblock On the performance degradation of cyber-physical systems under
  stealthy integrity attacks.
\newblock {\em IEEE Transactions on Automatic Control}, 61(9):2618--2624, 2016.

\bibitem{yilin4}
Yorie Nakahira and Yilin Mo.
\newblock Dynamic state estimation in the presence of compromised sensory data.
\newblock In {\em 54th Annual Conference on Decision and Control (CDC)}, pages
  5808--5813. IEEE, 2015.

\bibitem{pajic}
Miroslav Pajic, James Weimer, Nicola Bezzo, Paulo Tabuada, Oleg Sokolsky, Insup
  Lee, and George~J Pappas.
\newblock Robustness of attack-resilient state estimators.
\newblock In {\em ICCPS'14: ACM/IEEE 5th International Conference on
  Cyber-Physical Systems (with CPS Week 2014)}, pages 163--174, 2014.

\bibitem{fabio}
Fabio Pasqualetti, Florian Dorfler, and Francesco Bullo.
\newblock Attack detection and identification in cyber-physical systems.
\newblock {\em IEEE Transactions on Automatic Control}, 58(11):2715--2729,
  2013.

\bibitem{ricker93}
N~Lawrence Ricker.
\newblock Model predictive control of a continuous, nonlinear, two-phase
  reactor.
\newblock {\em Journal of Process Control}, 3(2):109--123, 1993.

\bibitem{sandberg}
Henrik Sandberg and Andr{\'e}~MH Teixeira.
\newblock From control system security indices to attack identifiability.
\newblock In {\em Science of Security for Cyber-Physical Systems Workshop
  (SOSCYPS)}, pages 1--6. IEEE, 2016.

\bibitem{paulo_persis2}
Danial Senejohnny, Pietro Tesi, and Claudio De~Persis.
\newblock A jamming-resilient algorithm for self-triggered network
  coordination.
\newblock {\em arXiv preprint arXiv:1603.02563}, 2016.

\bibitem{yasser_SMT}
Yasser Shoukry, Pierluigi Nuzzo, Alberto Puggelli, Alberto~L
  Sangiovanni-Vincentelli, Sanjit~A Seshia, and Paulo Tabuada.
\newblock Secure state estimation for cyber physical systems under sensor
  attacks: a satisfiability modulo theory approach.
\newblock {\em IEEE Transactions on Automatic Control}, 62(10):4917--4932,
  2017.

\bibitem{yasser1}
Yasser Shoukry and Paulo Tabuada.
\newblock Event-triggered state observers for sparse sensor noise/attacks.
\newblock {\em IEEE Transactions on Automatic Control}, 61(8):2079--2091, 2016.

\bibitem{showkatbakhsh3}
M.~Showkatbakhsh, Y.~Shoukry, H.~Chen R, S.~Diggavi, and P.~Tabuada.
\newblock An {SMT}-based approach to secure state estimation under sensor and
  actuator attacks.
\newblock In {\em Decision and Control (CDC), IEEE 56th Conference on}, pages
  7177--7182. IEEE, 2017.

\bibitem{showkatbakhsh2}
Mehrdad Showkatbakhsh, Paulo Tabuada, and Suhas Diggavi.
\newblock Secure system identification.
\newblock In {\em 54th Annual Allerton Conference on Communication, Control,
  and Computing}, pages 1137--1141. IEEE, 2016.

\bibitem{showkatbakhsh1}
Mehrdad Showkatbakhsh, Paulo Tabuada, and Suhas Diggavi.
\newblock System identification in the presence of adversarial outputs.
\newblock In {\em Decision and Control (CDC), IEEE 55th Conference on}, pages
  7177--7182. IEEE, 2016.

\bibitem{paulo_roysmith}
Roy~S Smith.
\newblock Covert misappropriation of networked control systems: Presenting a
  feedback structure.
\newblock {\em Control Systems Magazine, IEEE}, 35(1):82--92, 2015.

\bibitem{sundaram1}
Shreyas Sundaram, Miroslav Pajic, Christoforos~N Hadjicostis, Rahul Mangharam,
  and George~J Pappas.
\newblock The wireless control network: monitoring for malicious behavior.
\newblock In {\em 49th IEEE Conference on Decision and Control (CDC)}, pages
  5979--5984, 2010.

\bibitem{tiwari}
Ashish Tiwari, Bruno Dutertre, Dejan Jovanovi{\'c}, Thomas de~Candia, Patrick~D
  Lincoln, John Rushby, Dorsa Sadigh, and Sanjit Seshia.
\newblock Safety envelope for security.
\newblock In {\em ACM Proceedings of the 3rd international conference on High
  confidence networked systems}, pages 85--94, 2014.

\bibitem{yong}
Sze~Zheng Yong, Ming~Qing Foo, and Emilio Frazzoli.
\newblock Robust and resilient estimation for cyber-physical systems under
  adversarial attacks.
\newblock In {\em American Control Conference (ACC), 2016}, pages 308--315.
  IEEE, 2016.

\bibitem{partial}
T~Yoshikawa and S~Bhattacharyya.
\newblock Partial uniqueness: Observability and input identifiability.
\newblock {\em IEEE Transactions on Automatic Control}, 20(5):713--714, 1975.

\bibitem{paulo_sonia}
Minghui Zhu and Sonia Martinez.
\newblock On the performance analysis of resilient networked control systems
  under replay attacks.
\newblock {\em IEEE Transactions on Automatic Control}, 59(3):804--808, 2014.

\end{thebibliography}

\appendix

\begin{pf}[Lemma \ref{lemma:pre}]
We first prove the sufficiency part. For the sake of contradiction, suppose that the underlying system is not strongly observable but the property of Corollary \ref{cor:strong_observable0} is true. If the underlying system \eqref{sys} is not strongly observable, it means there exist two initial conditions, denoted by $x^{(1)}(0)$ and $x^{(2)}(0)$ possibly with different input sequences denoted by $\{ u^{(1)}(t) \}$ and $\{ u^{(2)}(t) \}$, respectively, that correspond to the same output sequence $\{y(t)\}$. The underlying system is linear, therefore the nonzero initial condition of $x^{(1)}(0) - x^{(2)}(0)$ with the input sequence $\{ u^{(1)}(t) - u^{(2)}(t) \}$ produces the zero output sequence which contradicts the property given in Corollary \ref{cor:strong_observable0}.
The necessity can be concluded using the similar argument. For the sake of contradiction let us assume this property does not hold, i.e., there exists a non zero initial state $x(0) \neq 0$ that corresponds to the zero output sequence. This contradicts the strong observability since the zero output sequence can be generated from both zero and $x(0) \neq 0$ as initial conditions under (possibly different) input sequences.
\end{pf}

\begin{pf}[Lemma \ref{lemma:sensors}]
We prove this lemma with contradiction. We show that if TEST($\Gamma_u^{\text{cert}}, \Gamma_y^{\text{temp}} \cup \{ i \}$) returns true for all $i \in \Gamma_y^{\text{SAT}} \setminus  \Gamma_y^{\text{temp}}$ then TEST($\Gamma_u^{\text{cert}}, \Gamma_y^{\text{SAT}}$) would also return true, which contradicts the assumption of the lemma. By applying the following lemma successively, the result follows directly.

\begin{lem}
Assume that the system $S$ is $(2r, 2s)$-sparse strongly observable. Pick any subset of inputs and outputs denoted by $\Gamma_u^{\text{cert}}$ and $\Gamma_y^{\text{temp}}$ with $|\Gamma_u^\text{cert}| \leq 2r$ and $|\Gamma_y^{\text{temp}}| \geq p-2s$. Then for any subsets of outputs denoted by $\Gamma_{y}^1$ and $\Gamma_{y}^2$, the first statement implies the second:
\begin{enumerate}
\item TEST($\Gamma_u^\text{cert}, \Gamma_y^{\text{temp}} \cup \Gamma_y^1 $) and TEST($\Gamma_u^\text{cert}, \Gamma_y^{\text{temp}} \cup \Gamma_y^2$) return true.
\item TEST($\Gamma_u^\text{cert}, \Gamma_y^{\text{temp}} \cup \Gamma_y^1 \cup \Gamma_y^2 $) returns true.
\end{enumerate}
\end{lem}

\begin{pf}
Without loss and generality we can assume $\Gamma_y^{1}$, $\Gamma_y^{2}$ and $\Gamma_y^\text{temp}$ are all disjoint sets. Since TEST($\Gamma_u^\text{cert}, \Gamma_y^{\text{temp}} \cup \Gamma_y^{i} $) returns true for $i \in \{1, 2\}$, therefore we have
\begin{align}
\begin{bmatrix} \bold{Y}\vert_{ \Gamma_y^{\text{temp}} } \\  \bold{Y}\vert_{\Gamma_y^1} \end{bmatrix} &= \begin{bmatrix} \mathcal{O}_{ \Gamma_y^{\text{temp}} } \\ \mathcal{O}_{\Gamma_y^1}  \end{bmatrix} { \hat{x}^{1} } + \begin{bmatrix} \mathcal{N}_{ \Gamma_u^{\text{cert}} \to \Gamma_y^{\text{temp}} } \\ \mathcal{N}_{ \Gamma_u^{\text{cert}} \to {\Gamma_y^1} } \end{bmatrix} \hat{\bold{U} }^{1}, \label{eq:proof:lemma:sensors:1} \\
\begin{bmatrix} \bold{Y}\vert_{ \Gamma_y^{\text{temp}} } \\  \bold{Y}\vert_{\Gamma_y^2} \end{bmatrix} &= \begin{bmatrix} \mathcal{O}_{ \Gamma_y^{\text{temp}} } \\ \mathcal{O}_{\Gamma_y^2}  \end{bmatrix} { \hat{x}^{2} } + \begin{bmatrix} \mathcal{N}_{ \Gamma_u^{\text{cert}} \to \Gamma_y^{\text{temp}} } \\ \mathcal{N}_{ \Gamma_u^{\text{cert}} \to {\Gamma_y^2} } \end{bmatrix} \hat{\bold{U} }^{2}, \label{eq:proof:lemma:sensors:2}
\end{align} where $\hat{x}^{1}, \hat{x}^2 \in \mathbb{R}^n$ are states that T-solver.check returns, $\hat{\bold{U}}^{1}, \hat{\bold{U}}^{2}$ are matrices with appropriate dimensions that satisfy TEST. 
Note that the underlying system is $(2r, 2s)$-sparse strongly observable, $| \Gamma_u^{\text{cert}} | \leq 2r$ and $| \Gamma_y^{\text{temp}} | \geq p-s$ therefore \mbox{$\hat{S} := \subsys{\Gamma_u^{\text{cert}} }{ \Gamma_y^{\text{temp}} }$} is strongly observable. One can reinterpret $(\hat{\bold{U} }^{1}, \bold{Y}\vert_{ \Gamma_y^{\text{temp}} })$ and $(\hat{\bold{U} }^{2}, \bold{Y}\vert_{ \Gamma_y^{\text{temp}} })$ as two (possibly different) valid trajectories of a strongly observable system $\hat{S}$ with identical output sequences. Strong observability implies that the state can be uniquely determined from the output with a delay bounded by $n$, therefore $\hat{x}^{1} = \hat{x}^2$. Furthermore, the equality of right hand sides of \eqref{eq:proof:lemma:sensors:1} and \eqref{eq:proof:lemma:sensors:2} implies that,
\begin{align}
\mathcal{N}_{ \Gamma_u^{\text{cert}} \to \Gamma_y^{\text{temp}} } (\hat{\bold{U} }^{2} - \hat{\bold{U} }^{1}) = {0}, \label{eq:proof:lemma:sensors:3}
\end{align} i.e., $\hat{\bold{U} }^{2} - \hat{\bold{U} }^{1}$ is a zero dynamic of $\hat{S}$. By $(2r, 2s)$-sparse strongly observablity of ${S}$, we conclude that $\hat{\bold{U} }^{2} - \hat{\bold{U} }^{1}$ is also a zero dynamic of $S$, and therefore,
\begin{align}
\mathcal{N}_{ \Gamma_u^{\text{cert}} \to \Gamma_y^{1} } (\hat{\bold{U} }^{2} - \hat{\bold{U} }^{1}) = {0}, \quad \mathcal{N}_{ \Gamma_u^{\text{cert}} \to \Gamma_y^{2} } (\hat{\bold{U} }^{2} - \hat{\bold{U} }^{1}) = {0} \label{eq:proof:lemma:sensors:4}.
\end{align} Putting \eqref{eq:proof:lemma:sensors:1}, \eqref{eq:proof:lemma:sensors:2} and \eqref{eq:proof:lemma:sensors:4} together with $\hat{x}^{1} = \hat{x}^2$, we conclude that:
\begin{align}
\begin{bmatrix} \bold{Y}\vert_{ \Gamma_y^{\text{temp}} } \\  \bold{Y}\vert_{\Gamma_y^1} \\  \bold{Y}\vert_{\Gamma_y^2}  \end{bmatrix} &= \begin{bmatrix} \mathcal{O}_{ \Gamma_y^{\text{temp}} } \\ \mathcal{O}_{\Gamma_y^1} \\ \mathcal{O}_{\Gamma_y^2}  \end{bmatrix} { \hat{x}^{1} } + \begin{bmatrix} \mathcal{N}_{ \Gamma_u^{\text{cert}} \to \Gamma_y^{\text{temp}} } \\ \mathcal{N}_{ \Gamma_u^{\text{cert}} \to {\Gamma_y^1} } \\ \mathcal{N}_{ \Gamma_u^{\text{cert}} \to {\Gamma_y^2} } \end{bmatrix} \hat{\bold{U} }^{1},
\end{align} i.e., TEST($\Gamma_u^\text{cert}, \Gamma_y^{\text{temp}} \cup \Gamma_y^{1} \cup  \Gamma_y^{2}$) returns false.

\end{pf}
\end{pf}

\begin{pf}[Proof Sketch of Lemma \ref{lemma:hueristics}]
Let us revisit the optimization \eqref{eq:test} inside the consistency check $\text{TEST}({ \Gamma_u }, {\Gamma_y} )$,
\begin{align}
 \min_{\hat{x}, \hat{\bold{U}}} \left\|
\bold{Y}\vert_{{\Gamma}_y} {-} 
\begin{bmatrix}
\mathcal{O}_{{\Gamma}_y}, \mathcal{N}_{{{\Gamma}_u} \to {{\Gamma}_y}}
\end{bmatrix}
\begin{bmatrix}
\hat{x} \\
\bold{\hat{U}}
\end{bmatrix} \right\| 
\end{align}
For a generic LTI system, the matrix $\begin{bmatrix} \mathcal{O}_{{\Gamma}_y}, \mathcal{N}_{{{\Gamma}_u} \to {{\Gamma}_y}} \end{bmatrix} \in \mathbb{R}^{ n |{ \Gamma }_y| \times n(1+| {\Gamma}_u| ) }$ is of full rank, where $n$ is the order of the LTI system. If $\begin{bmatrix} \mathcal{O}_{{\Gamma}_y}, \mathcal{N}_{{{\Gamma}_u} \to {{\Gamma}_y}} \end{bmatrix} \in \mathbb{R}^{ n |{ \Gamma }_y| \times n(1+| {\Gamma}_u| ) }$ is of full row rank, then $\text{TEST}({ \Gamma_u }, {\Gamma_y} )$ is satisfied irrespectively of the actual values of $\bold{Y}\vert_{{\Gamma}_y}$. Therefore in order to have a certificate constructed by inputs in $\overline{\Gamma}_u$ and outputs in $\Gamma_y$, $\begin{bmatrix} \mathcal{O}_{{\Gamma}_y}, \mathcal{N}_{{{\Gamma}_u} \to {{\Gamma}_y}} \end{bmatrix} \in \mathbb{R}^{ n |{ \Gamma }_y| \times n(1+| {\Gamma}_u| ) }$ should be a full column rank matrix, therefore
\begin{align}
n |{\Gamma}_y| \geq  n(1+| {\Gamma}_u| ). \label{eq:proof:heuristics:1} 
\end{align} 
The certificate consists of inputs in $\overline{{\Gamma}}_u$ and outputs in ${\Gamma}_y$, therefore the length of certificate is:
\begin{align}
| \overline{{\Gamma}}_u| + | {\Gamma}_y | &= m - | {\Gamma}_u | + | {\Gamma}_y |   \geq m+1.
\end{align}
\end{pf}

\end{document}